\documentclass[a4paper,12pt]{amsart}
\usepackage{hyperref,fancyhdr,amsmath,amscd,amsthm,amsfonts,latexsym,amssymb,mathrsfs}

\DeclareMathOperator{\trace}{trace}
\DeclareMathOperator{\grad}{grad}

\DeclareMathOperator{\dif}{d}

\renewcommand{\H}{\mathscr{H}}
\renewcommand{\o}{\mathcal{O}} 
\newcommand{\V}{\mathscr{V}}
\newcommand{\F}{\mathscr{F}} 
\newcommand{\K}{\mathscr{K}}  
\newcommand{\tG}{\mathscr{G}} 
 
\newcommand{\n}{\mathcal{N}}
 
\DeclareMathOperator{\Bh}{{\it B}^{\H}} 
\DeclareMathOperator{\UBh}{\mathcal{B}^{\H}} 
\DeclareMathOperator{\UB0h}{\mathcal{B}_0^{\H}} 
\DeclareMathOperator{\Bv}{{\it B}^{\V}}

\def \a{\alpha}

\def \b{\beta}

\def \g{\gamma}

\def \G{\Gamma}
\def \l{\lambda}
\def \O{\Omega}
\def \phi{\varphi} 
\def \p{\pi}
\def \r{\rho} 
\def \s{\sigma} 
\def \t{\tau}
\def \D{\Delta}
\def \R{\mathbb{R}}

\def \C{\mathbb{C}\,} 

\def\widecheckg{g^{\hspace*{-2.5pt}\vbox to 5pt{\hbox to
0pt{\LARGE$\check{}$}}}\hspace*{2pt}} 

\def\widecheckl{\lambda^{\hspace*{-3.5pt}\vbox to 8pt{\hbox to
0pt{\LARGE$\check{}$}}}\hspace*{2pt}}

\begin{document}

\title{Twistorial harmonic morphisms with one-dimensional fibres on self-dual four-manifolds}
\author{Radu Pantilie and John C.\ Wood} 
\thanks{The authors gratefully acknowledge that this work was 
done under E.P.S.R.C. grant number GR/N27897.} 
\email{Radu.Pantilie@imar.ro, J.C.Wood@leeds.ac.uk} 
\address{R.~Pantilie, Institutul de Matematic\v a ``Simion Stoilow'' al Academiei Rom\^ane, 
Calea Grivi\c tei nr. 21, 010702, C.P. 1-764, 014700, Bucure\c sti, Rom\^ania}
\address{J.C.~Wood, Department of Pure Mathematics, University of Leeds, 
Leeds LS2 9JT, England}

\newtheorem{thm}{Theorem}[section]
\newtheorem{lem}[thm]{Lemma}
\newtheorem{cor}[thm]{Corollary}
\newtheorem{prop}[thm]{Proposition}

\theoremstyle{definition}
\newtheorem{defn}[thm]{Definition}
\newtheorem{rem}[thm]{Remark}
\newtheorem{exm}[thm]{Example}

\numberwithin{equation}{section}

\maketitle
\thispagestyle{empty} 
\vspace{-0.5cm}
\section*{Abstract}
\begin{quote}
{\footnotesize  
We introduce a general notion of \emph{twistorial map} and classify twistorial
harmonic morphisms with one-dimensional fibres from self-dual four-manifolds. 
Such maps can be characterised as those which pull back Abelian monopoles to 
self-dual connections. In fact, the constructions involve solving a
\emph{generalised monopole equation}, and also the \emph{Beltrami fields equation} of 
hydrodynamics, and lead to constructions of self-dual metrics.}
\end{quote}   

\section*{Introduction}      

\indent
Harmonic morphisms between Riemannian manifolds are smooth maps which preserve
Laplace's equation. By the basic characterisation theorem (\,\cite{Fug}\,,
\cite{Ish}\,),   they are harmonic maps which are \emph{horizontally weakly
conformal} (see below).\\  \indent  Classification results for harmonic
morphisms with one-dimensional fibres can be found in  \cite{BaiWoo2}\,,
\cite{Bry}\,, \cite{Pan-thesis}\,, \cite{Pan-4to3} and \cite{PanWood-d}\,.  In
\cite{PanWood-d} it is proved that, from an Einstein manifold of  dimension
at least five, there are just \emph{two} types (\,\cite{BaiEel}, \cite{Bry}\,)
of  harmonic morphism with one-dimensional fibres. In dimension four, the
situation is different:  from an Einstein 4-manifold there are precisely
\emph{three} types of harmonic morphism with  one-dimensional fibres 
\cite{Pan-thesis}\,, \cite{Pan-4to3} (see also \cite{PanWoo-exm}\,), where the
first  two types are as before. It is significant that all these three types
of harmonic morphism are  \emph{twistorial maps} in the sense of Definition
\ref{defn:twistmap}, below. Moreover, by a result  of \cite{Woo-4d},
submersive harmonic morphisms from Einstein 4-manifolds to \mbox{Riemann}
surfaces  are twistorial maps.\\
\indent
We shall see that a submersion with (nowhere degenerate) one-dimensional fibres
from a four-dimensional (complex-)Riemannian manifold is twistorial if
and only if it is self-dual in the sense of \cite{Cal-sds}\,.\\
\indent
In this paper we classify
twistorial harmonic morphisms with one-dimensional fibres from real-analytic
Riemannian four-manifolds, finding precisely one more type which is related
to a metric construction of \cite{Cal-sds}\,.\\
\indent 
In Section \ref{section:facts}, we review some basic facts on
harmonic morphisms with  one-dimensional fibres. In Section
\ref{section:twiststr}, we introduce a notion of \emph{(almost) twistorial
structure}. Then, we recall (\,\cite{Pen}, \cite{AtHiSi}, 
\cite{Hit-complexmfds}, \cite{LeB-nullgeod}\,) some basic examples from 
twistor theory and show how they fit into our framework.
In Section \ref{section:twistmap}, we introduce the notion of
\emph{twistorial map} and
then, also in Sections \ref{section:nondegfibres} and \ref{section:4to3},  
we show that many examples and facts from twistor theory appear as natural
properties  of such twistorial maps. 
To make our definition of twistorial, we
must complexify our manifolds, so that  we shall often work with
complex-analytic maps between \emph{complex-Riemannian} manifolds (see
\cite{LeB-nullgeod}\,).\\
\indent 
The main result of Section \ref{section:4to3} is a reformulation
(Theorem \ref{thm:4to3}\,) of a  result of \cite{Cal-sds} which gives a nice
characterisation for twistorial maps with nowhere  degenerate fibres from a
four-dimensional oriented conformal manifold. We prove that the induced  Weyl
connection on the codomain coincides with the one hinted at in
\cite{Cal-sds}\,. Such twistorial maps pull back (Abelian) monopoles to
self-dual connections   (cf.\ \cite{CalPed}\,, \cite{Cal-sds}\,); by using 
Theorem \ref{thm:4to3}\,, we show that this property is \emph{equivalent} to the
property of being  twistorial. Then we show that twistorial harmonic
morphisms $\phi:(M^4,g)\to(N^3,h,D)$ are characterised  by the property that
there exist nontrivial monopoles on $(N^3,h,D)$ which are pulled back  to
flat connections.\\  
\indent  
In Section \ref{section:4to3harmorph} we prove
(Theorem \ref{thm:maint}\,) that, from a  real-analytic four-dimensional
Riemannian manifold, there are, up to conformal changes with basic factor, 
just \emph{four} types of twistorial harmonic morphism with one-dimensional
fibres, where the first  three types are as above with a slight extension
of type 3, and type 4 is new. The proof involves solving the 
\emph{monopole equation} \eqref{e:m} (cf.\ \cite{PanWoo-exm}\,).  Recall
\cite{Pan-thesis}\,,\,\cite{Pan-4to3}\,,\,\cite{PanWoo-exm} that harmonic
morphisms of type 3 are  determined by the Beltrami fields equation (see
\cite{KenPlu}\,).  Here, this equation appears once more: for a harmonic
morphism $\phi:(M^4,g)\to (N^3,h)$ of type 4\,,   the \emph{Lee form} $\a$
with respect to $h$ of the Weyl connection on $N^3$,  with respect to
which $\phi$ is twistorial, satisfies the Beltrami fields equation 
$\dif\!\a=\pm*\a$\, outside the zero set
of a function, up to a conformal change of $h$\,   (Proposition \ref{prop:c=pm1}\,).\\  
\indent 
Theorem \ref{thm:maint} together with a result of \cite{Cal-sds} 
(see Theorem \ref{thm:4to3}, below),  gives the classification of twistorial 
harmonic morphisms with one-dimensional fibres from a self-dual
four-dimensional manifold with real-analytic metric 
(Corollary \ref{cor:maintsd}).\\   
\indent 
In \cite{PanWoo-exm}, we gave a
new construction of Ricci-flat self-dual metrics based  on harmonic morphisms
of type 3. In Section \ref{section:newconstr}, we show that harmonic morphisms
are related to constructions of Einstein and self-dual metrics in \cite{JonTod} and \cite{Cal-sds}\,.\\
\indent 
We are grateful to D.M.J.~Calderbank and K.P.~Tod for useful comments on a 
preliminary  version of this paper.

\section{Some facts on harmonic morphisms with one-dimensional fibres} \label{section:facts}

\indent
In this section we present some basic facts on harmonic morphisms with one-dimensional fibres. 
See \cite{BaiWoo2} and \cite{Pan} for general accounts and \cite{Gudbib} for a list of
papers on harmonic morphisms. 

\begin{defn} 
A \emph{harmonic morphism} is a smooth map $\phi:(M^m,g)\to(N^n,h)$ between Riemannian manifolds which pulls back (locally defined) harmonic functions to harmonic functions, i.e., if 
$f:V \to \R$
is a harmonic function on an open subset of $N$ with $\phi^{-1}(V)$ nonempty,
then $f \circ \phi$ is a harmonic function on $\phi^{-1}(V)$. 
\end{defn} 

\begin{defn} \label{defn:horizontally conformal} 
A smooth map $\phi:(M^m,g)\to(N^n,h)$ between Riemannian manifolds is \emph{horizontally 
(weakly) conformal} if, at each point $x\in M$\,, \emph{either} $\dif\!\phi_x=0$\,, in which case 
$x$ is called a \emph{critical point} of $\phi$\,, \emph{or}\/   
$\dif\!\phi_x:T_xM\to T_{\phi(x)}N$ is surjective and its restriction to the 
horizontal space $\H_x=({\rm ker}\dif\!\phi_x)^{\perp}$ is a conformal (linear) map 
$(\H_x,g_x|_{\H_x})\to(T_{\phi(x)}N,h_{\phi(x)})$\,, in which case $x$ is called a
\emph{regular point} of $\phi$. Denote the conformality factor by $\l(x)$\,. 
The resulting function $\l$ is called the \emph{dilation} of $\phi$\,. The dilation 
is smooth outside the set of critical points and can be extended to a continuous function 
on $M^m$\,, with $\l^2$ smooth, by setting it equal to zero on the set of 
critical points.\\ 
\indent
A smooth map is called \emph{horizontally homothetic} if it is horizontally conformal 
with dilation constant along horizontal curves.\\ 
\indent
A \emph{homothetic foliation} is a foliation which is locally defined by horizontally 
homothetic submersions \cite{Pan}, \cite{Pan-thesis}\,.  
\end{defn} 

\begin{rem} 
A map $\phi:(M,g)\to(N,h)$ is horizontally weakly conformal if and only if,  
for each $x\in M$, the \emph{adjoint} $\dif\!\phi_x^*:(T_{\phi(x)}N,h_{\phi(x)})\to(T_xM,g_x)$ 
is a weakly conformal linear map (with image $\H_x$\,). 
This formulation shows that the condition of horizontal weak conformality is \emph{dual} to that 
of weak conformality (see also \cite[\S2.4]{BaiWoo2}\,).
\end{rem}  

\indent 
The basic characterisation result for harmonic morphisms is the following.  

\begin{thm}[\cite{Fug}, \cite{Ish}] \label{thm:FugIsh}  
A smooth map between Riemannian manifolds is a harmonic morphism if and only if it is a 
harmonic map which is horizontally weakly conformal. 
\end{thm} 

\indent
It follows that the set of regular points of a nonconstant harmonic morphism 
is an open dense subset of the domain \cite{Fug}\,.\\ 

\indent
The following two propositions give two of the four types of harmonic morphism with 
one-dimensional fibres that we shall meet (see Section 3 below). 

\begin{prop}[\cite{BaiEel}] \label{prop:BaiEel}  
Let $\phi:(M^{n+1},g)\to(N^n,h)$ be a nonconstant horizontally weakly conformal map 
between Riemannian manifolds of dimensions $n+1$ and $n$\,, respectively $(n\geq1)$\,. 
If $n=2$\,, then $\phi$ is a harmonic morphism if and only if its fibres are geodesic 
at regular points. If $n\neq2$\,, then any two of the following assertions imply the third:\\ 
\indent
{\rm (i)} $\phi$ is a harmonic morphism;\\ 
\indent
{\rm (ii)} the fibres of\/ $\phi$ are geodesic at regular points;\\ 
\indent
{\rm (iii)} $\phi$ is horizontally homothetic. 
\end{prop} 
\begin{proof} 
Let $\phi:(M^m,g)\to(N^n,h)$ be a horizontally weakly conformal submersion. 
Then, at regular points we have the following \emph{fundamental equation} (see, for example, \cite[\S4.5]{BaiWoo2}\,)
for the tension field $\t(\phi)\in\G(\phi^*(TN))$ of $\phi$\,: 
\begin{equation} \label{e:fundamentaleqn} 
\t(\phi)+\dif\!\phi(\trace(\Bv))=-(n-2)\dif\!\phi(\grad(\log\l)) 
\end{equation} 
where $\trace(\Bv)$ is the trace of the second fundamental form of $\V={\rm ker}\dif\!\phi$ 
and $\l$ is the dilation of $\phi$\,.\\ 
\indent 
The proposition is an immediate consequence of \eqref{e:fundamentaleqn}\,. 
\end{proof}  

\begin{rem} 
Proposition \ref{prop:BaiEel} is true for maps with higher-dimensional fibres 
after replacing `geodesic' with `minimal'. 
\end{rem} 

\indent
Recall that a foliation is called \emph{Riemannian} if it is locally 
defined by the fibres of Riemannian submersions; then we have the following result 
(see also \cite[\S12.3]{BaiWoo2}\,).

\begin{prop}[\cite{Bry}] \label{prop:Bry} 
Let $\phi:(M^{n+1},g)\to N^n$ be a submersion with connected one-dimensional fibres from a 
Riemannian manifold of dimension $n+1$ to a smooth manifold of dimension $n$ $(n\geq3)$\,. 
Suppose that the fibres of\/ $\phi$ form an orientable Riemannian foliation. Then the following 
assertions are equivalent:\\ 
\indent
{\rm (i)} There exists a Riemannian metric $h$ on $N^n$ with respect to which 
$\phi:(M^{n+1},g)\to(N^n,h)$ is a harmonic morphism;\\ 
\indent
{\rm (ii)} there exists a nowhere zero Killing vector field on $(M^{n+1},g)$ tangent 
to the fibres of\/ $\phi$\,. 
\end{prop} 

\indent
In general, harmonic morphisms with one-dimensional fibres have a local normal form 
as follows. 

\begin{thm}[\cite{Bry}] \label{thm:Bry} 
Let $(M^{n+1},N^n,S^1)$ be a principal bundle with projection $\phi:M^{n+1}\to N^n$,  
endowed with a principal connection $\H\subseteq TM$. Let $h$ be a Riemannian metric 
on $N^n$ and $\l$ a smooth positive function on $M^{n+1}$.\\ 
\indent
Define a Riemannian metric on $M^{n+1}$ by 
\begin{equation} \label{e:normalform} 
g=\l^{-2}\,\phi^*(h)+\l^{2n-4}\,\theta^2 
\end{equation} 
where $\theta$ is the connection form of\/ $\H$\,. Then  
$\phi:(M^{n+1},g)\to(N^n,h)$ is a harmonic morphism.\\ 
\indent
Conversely, any submersive harmonic morphism with one-dimensional fibres is \emph{locally} 
of this form, up to isometries. 
\end{thm} 

\indent
See \cite[\S12.2]{BaiWoo2}, \cite{Pan}, \cite{Pan-thesis} (and \cite{BaiWoo1} for the case 
$n=2$)  for a proof of Theorem \ref{thm:Bry} and a more explicit version of the converse.\\

\indent
The vector field $V$ on $M^{n+1}$ which is the vertical dual of\/ $\theta$ 
(i.e., $\dif\!\phi(V)=0$ and $\theta(V)=1$) is the infinitesimal generator of the local 
$S^1$-action; it is called \emph{the fundamental (vertical) vector 
field}; up to sign, it is characterised by the property that it is vertical and 
$g(V,V)=\l^{2n-4}$. 

\begin{rem} \label{rem:1.9} 
1) By Proposition \ref{prop:BaiEel}\,, any horizontally homothetic submersion with geodesic
fibres is a harmonic morphism. In the context of Theorem \ref{thm:Bry} this 
corresponds to the case when $\l$ is constant along horizontal curves; then if $\dif\!\l$ is nowhere 
zero, the horizontal distribution is integrable and totally umbilical.\\ 
\indent                              
2) The Killing vector field of Proposition \ref{prop:Bry} is equal to (a multiple of) the fundamental
vector field. In fact, from Proposition \ref{prop:Bry} and Theorem \ref{thm:Bry} 
(see \cite{Bry}, \cite[\S12.3]{BaiWoo2}, \cite{Pan}\,)
we deduce the following:  Let $V$ be a Killing vector field on  
$(M^{n+1},g)$\ ($n\neq2$)\  whose integral curves are the fibres of the submersion 
$\phi:(M^{n+1},g)\to N^n$; let $\widecheckg$ denote the unique metric on $N^n$ 
with respect to which $\phi:(M^{n+1},g)\to(N^n,\widecheckg\,)$ is a Riemannian submersion 
and $\widecheckl$ the positive smooth function on $N^n$ such that 
$\phi^*({\widecheckl}^{\;2n-4})=g(V,V)$\,; then $\phi:(M^{n+1},g)\to(N^n,{\widecheckl}^{\,-2}\,\widecheckg\,)$ is a harmonic morphism. 
In the context of Theorem \ref{thm:Bry} this corresponds to the case when $\l$ is constant 
along the fibres of\/ $\phi$\,. 
\end{rem} 

\indent
We end this section by recalling the following. 

\begin{defn}[\cite{Pan-thesis}, cf.\ \cite{Woo}] 
Let $\V$ be (the tangent bundle of) a foliation on the Riemannian manifold $(M^m,g)$\,.\\ 
\indent
We say that $\V$ \emph{produces harmonic morphisms} if it can be locally defined by submersive 
harmonic morphisms. 
\end{defn}

\section{Twistorial structures} \label{section:twiststr} 

\indent
In this section we shall work in the \emph{complex-analytic category}. Thus, all manifolds 
will be complex and all maps will be complex-analytic. Examples of such manifolds and maps 
can be obtained by complexifying real-analytic maps between real-analytic manifolds.
In particular, by a \emph{Riemannian manifold} we shall mean a \emph{complex-Riemannian manifold}
in the sense of \cite{LeB-nullgeod}\,.\\
\indent 
We define horizontally (weakly) conformal maps similarly to Definition \ref{defn:horizontally conformal}, 
but $\l$ may then be a complex number defined only up to sign with the \emph{square dilation} $\l^2$ globally defined.\\ 
\indent 
In what follows, it is convenient to work with the following definition.  

\begin{defn} \label{defn:atwiststr} 
Let $M$ be a (complex) manifold. By an \emph{almost twistorial structure (on the manifold $M$)} 
we shall mean a quadruple $(P, M, \p, \F)$ where $\p:P\to M$ is a proper surjective (complex-analytic)  
submersion and  $\F\subseteq TP$ is a distribution 
on $P$ such that $({\rm ker}\dif\!\p)\cap\F=\{0\}$\,; we call $\F$ the 
\emph{twistor distribution} of $(P, M, \p, \F)$\,.\\ 
\indent
We shall call the almost twistorial structure $(P, M, \p, \F)$ \emph{integrable} if its twistor 
distribution $\F$ is integrable. An integrable almost twistorial structure will be called a 
\emph{twistorial structure}. If $(P, M, \p, \F)$ is a twistorial structure then the leaf space of $\F$ 
is called the \emph{twistor space} of $(P, M, \p, \F)$\,.  
\end{defn} 

\begin{rem} \label{rem:family} 
1) Let $(P, M, \p, \F)$ be a twistorial structure. Assume that the foliation $\F$ is 
simple (that is, its leaf space $Z$ is a manifold and the projection $\pi_Z:P\to Z$ is 
a submersion) and that any of its leaves intersect each fibre of\/ $\pi$ at most once. 
Then$\{\p_Z(\pi^{-1}(x))\}_{x\in M}$ is an analytic family \cite{Kod} of compact 
submanifolds of\/ $Z$, which we shall call \emph{twistor submanifolds}, or 
\emph{twistor lines} when they are projective lines.\\
\indent
2) For all of the almost twistorial structures $(P,M,\p,\F)$\ $(\dim \F=k)$\ which will appear in this paper, 
the map $P\to G_k(TM)$\,, $p\mapsto \dif\!\p(\F_p)$\,, is an embedding of $P$ into the Grassmann bundle of 
$k$-dimensional planes on $TM$. Such almost twistorial structures can be obtained as follows.\\ 
\indent 
By a \emph{linear partial connection on a vector bundle $E\to M$,
over a distribution $\H$ on $M$} \cite{Bott-partial}\,, we mean a morphism $\nabla$ from 
the sheaf of sections 
of $E$ to the sheaf of sections of $E\otimes\H^*$ such that, if $U$ is an open set of $M$, 
then $\nabla:\G(U,E)\to\G(U,E\otimes\H^*)$ is a $\C$-linear map which satisfies  
$$\nabla(sf)=(\nabla s)f+s\otimes\dif\!f|_{\H}$$ 
for any function $f:U\to\C$ and section $s\in\G(U,E)$\,; in a 
similar way to the case of usual connections, any linear 
partial connection over $\H$ corresponds to a principal partial connection over $\H$ on the 
frame bundle of $E$ where 
by a \emph{principal partial connection, over $\H$\,, on the principal bundle $(P,M,G)$\,, with projection 
$\p:P\to M$}, we mean a $G$-invariant distribution $\K$ on $P$ such that $\K\cap{\rm ker}\dif\!\pi=\{0\}$ 
and $\dif\!\p(\K)=\H$\,.\\
\indent
Now let $\H\subseteq TM$ be an $n$-dimensional distribution endowed with a linear partial connection
$D$, over itself.
Suppose that $D$ corresponds to a principal partial connection, over $\H$, on a principal
subbundle $(Q,M,G)$, ($G\subseteq GL(n,\C)$\,), of the bundle of (complex) frames of $\H$\,. 
Let $F\subseteq G_k(\C^{\!n})$ be a submanifold which is invariant under the action of $G$ and 
let $P=Q\times_GF$ be the associated bundle. Clearly, $P\subseteq G_k(\H)$\,. Also, $D$ induces a 
connection $\mathcal{D}\subseteq TP$ on $\p:P\to M$. Then, for each $p\in P$ we 
define $\F_p\subseteq T_pP$ to be the horizontal lift, with respect to $\mathcal{D}$, 
at $p\in P$ of $p\subseteq \H_{\p(p)}$\,.\\ 
\indent 
Usually, $G\subseteq CO(n,\C)$ where $CO(n,\C)$ is the complex-conformal group in dimension $n$\,, 
so that $\H$ is endowed with a conformal structure. Then, if $D$ corresponds to a principal partial 
connection on the corresponding frame bundle $(Q,M,CO(n,\C))$\,, it is called a 
\emph{conformal partial connection}.\\
\indent
3) Let $(P,M,\p,\F)$, $(\dim M=m, \dim \F=k)$ be a twistorial structure where, as above, $P$ is a subbundle
of the Grassmann bundle $G_k(TM)$ and $\F$ is induced by some connection $D$ on $M$ which preserves $P$.
Suppose that $D$ is torsion-free and, for any $p\in P$, there exists a totally geodesic submanifold of $M$, of dimension $(m-k)$\,,
which passes through $\p(p)$ and which is transversal to $p$. Then each point of $M$ has an open
neighbourhood $U$ such that $\F|_{\p^{-1}(U)}$ is simple.        
\end{rem} 
 
\indent
Next, we give the basic examples of almost twistorial structures with which we shall work.
Recall that we are working in the complex-analytic category.  First, we consider structures
over $2$- and $3$-dimensional manifolds. 

\begin{exm}[C.R.~LeBrun \cite{LeB-nullgeod}] \label{exm:2t}
Let $M = M^2$ be a two-dimensional Riemannian manifold and let $\p:P\to M$ be the 
\emph{bundle of null directions} on $M^2$. Clearly, $P=\det(O(M))$ and hence $\p:P\to M$ is a 
$\mathbb{Z}_2$-covering. Furthermore, there exists a canonical one-dimensional foliation 
$\F$ on $P$ such that $\p$ locally maps leaves of $\F$ to (local) null geodesics on $M^2$. 
Hence any two-dimensional Riemannian (or conformal, if one prefers) manifold $M^2$
is canonically endowed with the twistorial structure $(P, M, \p, \F)$.\\ 
\indent
Conversely, any (almost) twistorial structure $(P,M,\p,\F)$ with $\dim M=2$\,, $\dim\F=1$ 
and $\p:P\to M$ a $\mathbb{Z}_2$-covering such that the map $P\to P(TM)$\,, 
$p\mapsto \dif\!\p(\F_p)$\,, is an embedding is induced by a conformal structure as above.\\  
\indent 
If the Riemannian manifold $M^2$ is orientable (equivalently, if the $\mathbb{Z}_2$-covering 
$\p:P\to M$ is trivial) then $P=M_+\sqcup M_-$ where 
$M_+$ and $M_-$ are copies of $M$\,. Hence, there exist two foliations by null geodesics 
$\F_+$ and $\F_-$, on $M$, which are the projections of $\F$ restricted to $M_+$ and $M_-$, respectively. 
Therefore, the twistor space $Z=Z(M)$ of the canonical twistorial structure 
of an oriented two-dimensional Riemannian manifold $M^2$ (more generally, of a two-dimensional manifold 
endowed with an oriented conformal structure) is the \emph{space of null geodesics}  
\cite{LeB-nullgeod} of $M^2$. Then, locally, $Z$ is 
the disjoint union of two curves $Z=C_+\sqcup C_-$ and $M=C_+\times C_-$. 
Also, note that the (complex-analytic) almost Hermitian structures $J_{\pm}$ on $M$ defined by 
$J_{\pm}(X)=\pm{\rm i}X$ for $X\in\F_{\pm}$ are integrable. 
\end{exm} 

\indent
Let $(M,g)$ be a Riemannian manifold, $\dim M=m$\,. A \emph{degenerate hyperplane} $H\subseteq TM$ is a subspace of codimension one such that $g|_H$ is degenerate (equivalently, $H$ is the orthogonal complement of a null 
vector which is thus contained in $H$) (cf.\ \cite{LeB-nullgeod}\,); if $\dim M=3$ 
we shall say \emph{degenerate plane}. For fixed $x\in M$, the space of degenerate 
hyperplanes in $T_xM$ can be identified with the nonsingular quadric $Q_{m-2}$ 
in $P(T_xM)$.  

\begin{exm}[N.J.~Hitchin \cite{Hit-complexmfds}] \label{exm:3t} 
Let $M^3$ be a three-dimensional Riemannian manifold. Let $D$ be a Weyl connection  
(that is, a torsion-free conformal connection) on $M$. 
Let $\p:P\to M$ be the \emph{bundle of degenerate planes} on $M^3$. Then $P=CO(M)\times_{\r}\C\!P^1$ 
where $CO(M)$ is bundle of conformal frames of $M^3$ and $\r$ is the action of $CO(3,\C)$ on 
$\C\!P^1=\{\,p\,|\,p\subseteq \C^3 \,{\rm degenerate\:plane}\,\}$ 
induced by the canonical action of $CO(3,\C)$ on $\C^3$\,. Thus $D$ induces a connection on $\p:P\to M$ 
and, for each $p\in P$, we define $\F_p\subseteq T_pP$ to be the horizontal lift at 
$p\in P$ of $p\subseteq T_{\p(p)}M$\,. Obviously, the almost twistorial structure 
$(P, M, \p, \F)$ depends only on $D$ and on the conformal class $c$ of the metric of $M^3$. 
 
\begin{thm}[\cite{Hit-complexmfds}] \label{thm:3t} 
The twistor distribution $\F$ is integrable if and only if $D$ is Einstein--Weyl.\\
\indent
Furthermore, if $D$ is Einstein--Weyl then (locally) the leaf space $Z$ of\/ $\F$ contains a 
locally complete analytic family of projective lines each of which has normal 
bundle $\o(2)$\,. Conversely, any surface $Z$ containing a projective line with 
normal bundle $\o(2)$ is (locally) the twistor space of a $3$-dimensional Riemannian 
manifold $M^3$ endowed with an Einstein--Weyl connection $D$. The conformal structure 
of\/ $M^3$ and the Einstein--Weyl connection $D$ are uniquely determined.   
\end{thm} 

\indent
Let $(P, M, \p, \F)$ be the twistorial structure corresponding, as above, to $(M^3,D)$ where $D$ 
is an Einstein--Weyl connection on $M^3$. Then, as just explained, $\F$ is integrable and the 
locally complete analytic family \cite{Kod} of projective lines on $Z$ of 
Theorem \ref{thm:3t} 
appears as in Remark \ref{rem:family}(1) (see, also, Remark \ref{rem:family}(3)\,): 
each projective line represents all the (local) degenerate surfaces through a given point 
which are totally-geodesic with respect to $D$\,. 
The fact that the normal bundle of any twistor line in $Z$ is $\o(2)$ can be proved as follows. 
Let $t=\C\!P^1$ be a fibre of $\pi$\,. Then the normal bundle of $\pi_Z(t)=\C\!P^1$ in $Z$ 
is (isomorphic to) $(t\times \C^3)\big/(\F|_t)$\,. Now, $\F|_t$ is the 
restriction to $t=Q_1\subseteq\C\!P^2=G_2(\C^3)$ of the tautological plane bundle $E$ over 
$G_2(\C^3)$ where $Q_1$ is the one-dimensional nonsingular quadric and $G_2(\C^3)$ is the 
Grassmann manifold of planes in $\C^3$. As $(G_2(\C^3)\times\C^3)\big/E=H$ where $H\to\C\!P^2$ is 
the hyperplane bundle and the embedding $t=\C\!P^1=Q_1\hookrightarrow\C\!P^2$ has degree two 
we obtain $(t\times \C^3)\big/(\F|_t)=\o(2)$\,.\\ 
\indent
Note that, locally, any leaf of $\F$ is mapped by $\p$ to a degenerate surface in $M^3$ which is 
totally geodesic with respect to $D$, and so $Z$ is, locally, the space of degenerate  surfaces on 
$M^3$ which are totally geodesic with respect to $D$\,. 
\end{exm}

\indent 
Finally, we discuss an important example of almost twistorial structures over a four-dimensional manifold.

\begin{exm}[R.~Penrose \cite{Pen}; M.F.~Atiyah, N.J.~Hitchin, I.M.~Singer \cite{AtHiSi}] \label{exm:4t} 

\indent   
(i) Let $M^4$ be a four-dimensional Riemannian manifold. A plane $H\subseteq TM$ on $(M^4,g)$ 
is called \emph{null} if $g|_H=0$\,. Let $\p:P\to M$ be the \emph{bundle of null planes} 
on $M^4$. Then the Levi-Civita connection of $M^4$ induces a connection on $\p:P\to M$  
and, for each $p\in P$, we define $\F_p\subseteq T_pP$ to be the horizontal lift at  
$p\in P$ of $p\subseteq T_{\p(p)}M$\,. As $\F$ is conformally invariant, the almost twistorial structure 
$(P, M, \p, \F)$ canonically associated to $M^4$ is conformally invariant. It is well known that $\F$ 
is integrable if and only if $M^4$ is conformally flat (a consequence of 
Theorem \ref{thm:PenAtHiSi-4t}, below).\\           
\indent
(ii) We have that $(M^4,g)$ is orientable if and only if $P=P_+\sqcup P_-$ is the disjoint union of two $\C\!P^1$-bundles {over $M^4$}. In this case,
with respect to a choice of orientation on $M^4$, $P_+\to M$ is the \emph{bundle of
self-dual planes} on $M^4$ and $P_-\to M$ is the \emph{bundle of anti-self-dual planes} on $M^4$. 
(Self-dual and anti-self-dual planes are sometimes called \emph{$\a$-planes} and 
\emph{$\b$-planes\/}, see \cite{WarWel}, \cite{MasWoo}\,.) 
Then, with $\p_{\pm}=\p|_{P_{\pm}}$ and $\F_{\pm}=\F|_{P_{\pm}}$\,, $(P_{\pm}, M, \p_{\pm}, \F_{\pm})$ 
are almost twistorial structures on $M^4$. 
These almost twistorial structures are conformally invariant and so are well defined on any 
four-dimensional manifold endowed with an oriented conformal structure.  

\begin{thm}[\cite{Pen}, \cite{AtHiSi}] \label{thm:PenAtHiSi-4t} 
The twistor distribution $\F_-$ is integrable if and only if $M^4$ is self-dual. 
(Similarly, $\F_+$ is integrable if and only if $M^4$ is anti-self-dual.)\\ 
\indent
Furthermore, if $M^4$ is self-dual then (locally) the leaf space $Z$ of\/ $\F_-$ contains a 
locally complete analytic family of projective lines each of which has normal 
bundle $\o(1)\oplus\o(1)$\,. Conversely, any three-dimensional manifold $Z$ containing a projective 
line with normal bundle $\o(1)\oplus\o(1)$ is (locally) the twistor space of a self-dual Riemannian 
$4$-manifold which is uniquely determined, up to a conformal change of the metric. 
\end{thm} 

\indent
As the oriented $\mathbb{Z}_2$-covering of a \emph{non}orientable Riemannian manifold is canonically 
endowed with an orientation-reversing isometry, the oriented $\mathbb{Z}_2$-covering of a 
nonorientable Riemannian four-manifold $M^4$ is (anti-)self-dual if and 
only if $M^4$ is conformally flat.\\ 
\indent
Let $(P_-, M, \p_-, \F_-)$ be the twistorial structure corresponding, as above, to the four-dimensional 
self-dual Riemannnian manifold $M^4$. Then, by Theorem \ref{thm:PenAtHiSi-4t}, $\F_-$ is integrable  
and the locally complete analytic family of projective lines \cite{Kod} on $Z$ appears again as in 
Remark \ref{rem:family}(1) (see, also, Remark \ref{rem:family}(3)\,): each projective 
line represents all the (local) anti-self-dual surfaces through a given point.  
The fact that the normal bundle of any twistor line in $Z$ is $\o(1)\oplus\o(1)$ can be proved as follows. 
Let $t=\C\!P^1$ be a fibre of $\pi_-$\,. Then the normal bundle of $\pi_Z(t)=\C\!P^1$ in $Z$ 
is (isomorphic to) $(t\times \C^4)\big/(\F_-|_t)$\,. Now, $\F_-|_t$ is the 
restriction to $t=\C\!P^1\subseteq G_2(\C^4)$ of the tautological plane bundle $E$ over 
$G_2(\C^4)$. As any anti-self-dual $p$ plane on $\C^4=\C^2\otimes\C^2$ is given by 
$p=\{\,v\otimes v'\,|\,v\in\C^2\,\}$ for some fixed $v'\in\C^2$ we have that $E|_t=\o(-1)\oplus\o(-1)$ 
and hence the normal bundle of $\p_Z(t)$ in $Z$ is $\o(1)\oplus\o(1)$\,.\\ 
\indent
Note that, locally, any leaf of $\F_-$ is mapped by $\p_-$ to an anti-self-dual surface in $M^4$  
so $Z$ is, locally, the space of anti-self-dual surfaces in $M^4$. 
\end{exm}

\begin{rem}
The twistorial structures of Example \ref{exm:3t} are \emph{reductions} of the twistorial structures
of Example \ref{exm:4t} \cite{Hit-complexmfds} (see \cite{Cal-sds}\,).
\end{rem}

\section{Twistorial maps} \label{section:twistmap} 

\indent
In this section, unless otherwise stated, we again work in the complex-analytic category, thus 
all the manifolds are complex and all the maps are complex-analytic.\\ 
\indent 
In what follows, it is convenient to work with the following definitions. We start 
with an important special case.   

\begin{defn} \label{defn:tm} 
Let $\s=(P, M, \p_P, \F)$ and $\t=(Q, N, \p_Q, \tG)$ be almost twistorial structures 
over $M$ and $N$, respectively. Let $\phi:M\to N$ be a map. Suppose that there exists a map
$\Phi:P\to Q$ which covers $\phi$ (that is, $\p_Q\circ\Phi=\phi\circ\p_P$).\\
\indent
Then $\phi:(M,\s)\to(N,\t)$ will be called a \emph{twistorial map (with respect to $\Phi$)}
if $\dif\!\Phi(\F_p)\subseteq\tG_{\Phi(p)}$ for all $p\in P$. 
\end{defn} 

\indent
The following lemma allows us to construct `substructures' of an almost twistorial structure; its proof 
is omitted. 

\begin{lem} \label{lem:tm} 
Let $\s=(P, M, \p, \F)$ be an almost twistorial structure and let $P'\subseteq P$ be a closed 
submanifold such that $\p(P')=M$. Suppose that $\F_p\subseteq T_pP'$ for all $p\in P'$ and 
$\dim(T_pP'\cap{\rm ker}\dif\!\p_p)$ does not depend on $p\in P'$.\\
\indent
Then $\s'=(P', M, \p|_{P'}, \F|_{P'})$ is an almost twistorial structure. 
\end{lem}  

\indent
Next, we generalize the definition of twistorial map by allowing $\Phi$ to be a map between 
such substructures.  

\begin{defn} \label{defn:twistmap} 
Let $\s=(P, M, \p_P, \F)$ and $\t=(Q, N, \p_Q, \tG)$ be almost twistorial structures 
over $M$ and $N$, respectively, and let $\phi:M\to N$ be a map. Let $P'\subseteq P$ be as in 
Lemma \ref{lem:tm} and similarly for $Q'\subseteq Q$ so that $\s'=(P', M, \p_P|_{P'}, TP'\cap\F)$
and $\t'=(Q', N, \p_Q|_{Q'}, TQ'\cap\tG)$ are almost twistorial structures. Suppose that $\phi$ is
covered by a map $\Phi:P'\to Q'$. Then the map $\phi:(M,\s)\to(N,\t)$ will be called a
\emph{twistorial map (with respect to $\Phi$)}
if $\dif\!\Phi(\F_p)\subseteq \tG_{\Phi(p)}$ for all $p\in P'$.\\
\indent
Often, our choice of $P'$ and $Q'$ will depend on the map $\phi$\,.\\
\indent 
If the distributions $\F|_{P'}$ and $\tG|_{Q'}$ are integrable with leaf spaces 
$Z(M,\phi)$ and $Z(N,\phi)$, respectively, then we have an induced local map 
$Z(\phi):Z(M, \phi)\to Z(N, \phi)$ which we shall call the \emph{twistorial representation} of $\phi$. 
\end{defn} 

\begin{rem} \label{rem:tm}
1) Any real-analytic (Riemannian) manifold has a germ-unique complexification to a complex(-Riemannian) manifold \cite{LeB-nullgeod}.  We can \emph{complexify}
any real-analytic map between real-analytic manifolds, i.e., extend it to (the germ of) a complex-analytic map between
complexifications of those manifolds.
A map between real-analytic manifolds will be called twistorial if its complexification is twistorial.\\
\indent
2) Note that Definition \ref{defn:tm} is the case when $P'=P$ and $Q'=Q$\,.\\
\indent 
3) In all our examples, $P$ is embedded as a subbundle
of some Grassmann bundle $G_k(TM)$\,, as in Remark \ref{rem:family}(2)\,.
Also, $\Phi$ will be naturally induced by $\dif\!\phi$ so that we shall write $\Phi=\dif\!\phi|_{P'}$\,.\\ 
\indent 
4) In Definition \ref{defn:twistmap}, if the distributions $\F$ and $\tG$ are integrable
and $Z(M)$, $Z(M,\phi)$, $Z(N)$, $Z(N,\phi)$ are the leaf spaces of $\F$\,, $\F|_{P'}$\,, $\tG$\,, 
$\tG|_{Q'}$\,, respectively, then we have induced local maps $Z(M, \phi)\to Z(M)$, $Z(N, \phi)\to Z(N)$\,.\\ 
\indent
5) Most of the twistorial maps which will appear in this paper are submersive maps.
We could also consider \emph{twistorial foliations}, that is, foliations which are locally defined by 
submersive twistorial maps; most of the facts that follow can be easily reformulated in terms 
of such foliations.\\ 
\indent
6) {\bf (Compositions)} If $M_j$ \ $(j=1, 2, 3)$\ are endowed with almost twistorial structures
$\t_j=(P_j,M_j,\p_j,\F_j)$ \ $(j=1, 2, 3)	$\ and $\phi_1:(M_1, \t_1)\to(M_2, \t_2)$, 
$\phi_2:(M_2, \t_2)\to(M_3, \t_3)$ are twistorial maps with respect to $\Phi_1:P'_1\to Q'_1$ and
$\Phi_2:P'_2\to Q'_2$\,, such that $\Phi_1(P'_1)\subseteq P'_2$\,,
then $\phi_2\circ\phi_1:(M,\t_1)\to(M,\t_3)$ is a twistorial map with respect to 
$\Phi_2\circ\Phi_1$\,. 
\end{rem} 

\indent
The following simple proposition will be useful later on. 

\begin{prop} \label{prop:Grassmann} 
Let $M$ be a manifold endowed with a linear connection $D$ and let 
$\H\subseteq T(G_k(TM))$ be the connection induced by $D$ on the Grassmann bundle 
$\p:G_k(TM)\to M$ of\/ $k$-dimensional planes on $M$ ($k\leq\dim M$). Let $\F$ be the $k$-dimensional 
distribution on $G_k(TM)$ defined by setting $\F_p$ equal to the horizontal lift at 
$p\in G_k(TM)$ of\/ $p\subseteq T_{\p(p)}M$.\\ 
\indent 
Then, for a distribution $s:M\to G_k(TM)$ on $M$ the following assertions are equivalent:\\ 
\indent
{\rm (i)} $\F_{s(x)}\subseteq T_{s(x)}(s(M))$ for all $x\in M$;\\ 
\indent
{\rm (ii)} for any curve $c$ tangent to the distribution $s$ on $M$  
we have that $s\circ c$ is a parallel section of\/ $G_k(TM)$;\\ 
\indent
{\rm (iii)} for any vector fields $X,Y\in\G(s)$ we have $D_XY\in\G(s)$\,.\\ 
\indent
Furthermore, if $D$ is torsion free then, the following assertion can be added to this list:\\ 
\indent
{\rm (iv)} the distribution $s$ is integrable and its leaves are totally geodesic 
with respect to $D$\,. 
\end{prop} 
\begin{proof} 
Assertion (i) is equivalent to the fact that for any $X\in s(x)$ we have  
$\dif\!s(X)\in\F_{s(x)}$\,; this is clearly equivalent to assertion (ii).\\ 
\indent
Assertion (ii) is equivalent to the fact that the distribution $s$ is invariant under 
parallel transport along curves tangent to $s$\,. Then, the equivalence of (ii) and (iii) 
follows from the fact that, if $U\subseteq\C$ is a domain in $\C$, $H$ is a Lie subgroup 
of the Lie group $G$ and $(a_z)_{z\in U}$ is a curve in $G$ such that $a_{z_0}\in H$ for some $z_0\in U$ 
and $\theta(\dif\!a_z/\dif\!z)\in L(H)$ for all $z\in U$, where 
$\theta\in\G(L(G)\otimes T^*G)$ is the canonical form of\/ $G$, then $a_z\in H$ for all $z\in U$.\\ 
\indent
If $D$ is torsion free, the equivalence (iii)$\iff$(iv) is trivial. 
\end{proof} 

\begin{exm} \label{exm:2to4} 
We interpret the twistor lift of J.~Eells and S.~Salamon \cite{EelSal} in our framework. 
Let $M^2$ and $N^4$ be manifolds of dimension two and four, respectively, endowed with oriented conformal 
structures and let $\phi:M^2\to N^4$ be an injective conformal immersion.\\ 
\indent 
Endow $M^2$ and $N^4$ with the almost twistorial structures $(P,M,\p,\F)$ and $(P_-,N,\p_-,\F_-)$ of 
Examples \ref{exm:2t} and \ref{exm:4t}\,, respectively.\\ 
\indent 
As $\phi$ maps null directions on $M^2$ to null directions on $N^4$, we can define $\Phi_-:P\to P_-$ 
by $\Phi_-(p)$ is the anti-self-dual plane containing $\dif\!\phi(p)$\,, for each null direction $p$ on $M^2$. 
Then, $\phi$ is twistorial with respect to $\Phi_-$ if and only if, for any null geodesic $\g$ on $M^2$,  
$\Phi_-\circ\g$ is parallel along $\phi\circ\g$\,. If $N^4$ is self-dual then $\phi$ is twistorial 
with respect to $\Phi_-$ if and only if, for any null geodesic $\g$ on $M^2$ there exists a 
(necessarily unique) anti-self-dual surface $S_{\g}\subseteq N^4$ such that $\phi(\g)\subseteq S_{\g}$\,. 
Then, the map $\g\mapsto S_{\g}$ is the twistorial representation of\/ $\phi$ whose image is a pair of 
(local) curves in $Z(N)$\,.\\ 
\indent 
A nowhere degenerate surface in $N^4$ will be called a \emph{$(-)$twistorial surface} if 
the corresponding inclusion map is twistorial, in the above sense.\\
\indent
Similarly, we define \emph{$(+)$twistorial surfaces} by using instead the almost twistorial structure
$(P_+,N,\p_+,\F_+)$ of Example \ref{exm:4t}(ii)\,.
\emph{A nowhere degenerate surface on
an oriented four-dimensional conformal manifold is totally umbilical if and only if it is both 
$(+)$twistorial and $(-)$twistorial} \cite{EelSal}\,.    
\end{exm}

\section{Twistorial maps to surfaces} \label{section:nondegfibres}

\indent 
We shall now describe the twistorial maps which we shall need, namely those with nowhere degenerate 
fibres.

\begin{exm}[cf.\ \cite{BaiWoo-Bernstein}, \cite{BaiWoo-spfo}, \cite{CalPed}] \label{exm:3to2}  
Let $\phi:M^3\to N^2$ be a surjective submersion whose fibres are nowhere null from a 
Riemannian manifold onto an oriented Riemannian manifold. Let $D$ be 
a Weyl connection on $M^3$. We consider $M^3$ to be endowed with the almost twistorial 
structure $(P,M,\p,\F)$ of Example \ref{exm:3t}\,, and $N^2$ with the twistorial structure of 
Example \ref{exm:2t}(i) with $Q\to N$ the bundle of null directions on $N^2$.\\ 
\indent
At each $x\in M$ there are precisely two horizontal null directions $h_+(x)$ and $h_-(x)$\,. 
Thus to $\phi$ correspond the two disjoint embeddings $M_{\pm}=h_{\pm}(M)\hookrightarrow P$
given by $x\mapsto h_{\pm}(x)^{\perp}={\rm Span}(h_{\pm}(x), {\rm ker}\dif\!\phi_x)$\,.\\ 
\indent
As $\dim P=4$, $\dim\F=2$ and $\dim M_{\pm}=3$\,, at each degenerate plane $p\in M_{\pm}$, 
$d_{\pm}(p)=\dim(\F_p\cap T_pM_{\pm})\in\{1, 2\}$\,.\\ 
\indent
By Proposition \ref{prop:Grassmann}, $d_{\pm}=2$ on $M_{\pm}$ if and only if 
$\phi:M^3\to N^2$ is a horizontally conformal submersion with fibres which are
geodesic with respect to $D$. Therefore, $\phi$ is a twistorial map with $P'=M_+\sqcup M_-$ 
(and $\Phi=\dif\!\phi|_{P'}:P'\to Q$)  
if and only if it is a horizontally conformal submersion with geodesic fibres 
(with respect to $D$)\,. If $D$ is Einstein--Weyl, then 
$Z(M,\phi)=C_+\sqcup C_-$ is the leaf space of the foliation induced by $\F$ on $M_+\sqcup M_-$  
and the map $Z(M,\phi)\to Z(M)$ is simply the inclusion $C_+\sqcup C_-\hookrightarrow Z(M)$\,. 
Obviously $C_+$ and $C_-$ are transversal to the twistor lines in $Z(M)$\,.
Also, $Z(N,\phi)=Z(N)$ is the space of null geodesics of $N^2$, which, as a consequence of the 
horizontal conformality of $\phi$ can be canonically identified with $C_+\sqcup C_-$\,, then 
the twistorial representation of $\phi$ is simply the identification map 
$C_+\sqcup C_- \to Z(N)$. To retrieve $\phi$ from its twistorial representation, let $x\in M^3$;  
this corresponds to the twistor line $t_x\subseteq Z(M)$. Then, locally, $t_x$ meets $C_+$ and 
$C_-$ in two points which, under the identification $C_+\sqcup C_-=Z(N)$, correspond to two null 
geodesics on $N^2$ whose intersection is precisely $\phi(x)$.\\ 
\indent                                       
If $d_{\pm}=1$ on $M_{\pm}$ then $\F$ induces a one-dimensional foliation on $M_+\cup M_-$. 
For example, if $\phi:M^3\to N^2$ is a horizontally conformal submersion (with nowhere null 
fibres) whose fibres are nowhere geodesics, then $d_{\pm}=1$ on $M_{\pm}$ and the leaves
of the induced foliation on $M_+\sqcup M_-$ project onto the horizontal null geodesics. 
\emph{Therefore, a horizontally conformal submersion $\phi:M^3\to N^2$ is a twistorial map 
if and only if its fibres are geodesics (with respect to $D$\,).}
For such a horizontally conformal submersion, if $D$ is Einstein--Weyl, then, at least locally, 
any totally geodesic (with respect to $D$) degenerate surface contains precisely one horizontal null geodesic 
and so we have two \emph{local} sections $Z_{\pm}\hookrightarrow\n(M)$ of the canonical projection 
\cite{Hit-complexmfds} $\n(M)\to Z(M)$ where $\n(M)$ is the space of null geodesics of $M^3$ so that 
$Z_+\sqcup Z_-$ is the space of null geodesics of $M^3$ which are horizontal with respect to $\phi$\,. 
If $\phi$ is twistorial (with respect to $D$), then the fibres of $\n(M)\to Z(M)$ induce a 
one-dimensional foliation on $Z_+\sqcup Z_-$ whose leaf space is $C_+\sqcup C_-$\,.   
\end{exm} 

\indent
The following result (cf.\ Theorem \ref{thm:4to2harmorph} and 
Theorem \ref{thm:4to3harmorph}, 
below), which is a consequence of Proposition \ref{prop:BaiEel} and Example \ref{exm:3to2}, 
is a rephrasing of the starting point of the classification results of \cite{BaiWoo-Bernstein}, 
\cite{BaiWoo-spfo} (see also \cite[Chapters 1 and 6]{BaiWoo2}), there given for the smooth category. 

\begin{prop} \label{prop:3to2harmorph} 
Let $\phi:M^3\to N^2$ be the complexification of a real-analytic submersive map between 
Riemannian manifolds. Equip $M^3$ with the Levi-Civita connection and endow it with the almost 
twistorial structure of Example \ref{exm:3t}\,; equip $N^2$ with the twistorial structure of 
Example \ref{exm:2t}\,.\\ 
\indent 
Then, $\phi$ is a twistorial map, in the sense of Example \ref{exm:3to2}\,, if and only if 
it is the complexification of a harmonic morphism. \qed
\end{prop}  

\begin{defn} \label{defn:(-)twist}
Let $M^4$ be an oriented four-dimensional Riemannian manifold considered with the almost twistorial 
structure $\t_-=(P_-,M,\p_-,\F_-)$ of Example \ref{exm:4t}(ii)\,, with $\p_-:P_-\to M$ the bundle of 
anti-self-dual planes. Let $\phi:M^4\to N^n$ be a surjective submersive map with $n \ (=\dim N)=2$ or $3$\,. 
Let $\t =(Q,N,\p_,\,\tG)$ be the almost twistorial structure of Example \ref{exm:2t} (if $n=2$), or of 
Example \ref{exm:3t} for some Weyl connection on $N$ (if $n=3$). 

We shall say that $\phi$ is \emph{$(-)$twistorial (with respect to $\Phi$)} if $\phi:(M,\t_-)\to(N,\t)$ is a twistorial map with respect to $\Phi=\dif\!\phi|_{P'}:P'\to Q$\ for some $P' \subseteq P_-$ as 
in Lemma \ref{lem:tm}. 

Similarly, we define what is meant by $\phi$ is \emph{$(+)$twistorial} by using instead the almost twistorial structure $(P_+,N,\p_+,\F_+)$ of Example \ref{exm:4t}(ii)\,.
We shall write \emph{$(\pm)$twistorial} to mean \emph{$(+)$twistorial or $(-)$twistorial}.

\end{defn}  

\begin{exm}[cf.\ \cite{Woo-4d}] \label{exm:4to2}
Let $\phi:M^4\to N^2$ be a surjective submersion between oriented Riemannian manifolds, 
whose fibres are nowhere degenerate. We consider $M^4$ endowed with the almost twistorial 
structure $(P_-,M,\p_-,\F_-)$ of Examples \ref{exm:4t}(ii) with $P=P_-$ the bundle of
anti-self-dual planes on $M^4$, and $N^2$ endowed with the twistorial structure $(Q,N,\p,\tG)$ of
Examples \ref{exm:2t}\,, with $Q\to N$ the bundle of null directions on $N^2$.\\ 
\indent
At each point $x\in M^4$ there are precisely two horizontal null directions $h_+(x)$, $h_-(x)$ 
and two vertical null directions $v_+(x)$, $v_-(x)$. 
Let $p_{+}(x)$ (respectively, $p_{-}(x)$\,) be the (null) plane spanned by 
$h_{+}(x)$ and $v_{+}(x)$ (respectively, $h_{-}(x)$\, and $v_{-}(x)$\,). 
Then, from the fact that through each null direction there passes a unique (anti-)self-dual plane,  
it easily follows that either both $p_+(x)$ and $p_-(x)$ are self-dual or both are anti-self-dual. 
We assume that we have chosen the orientations such that both $p_+(x)$ and $p_-(x)$ are 
anti-self-dual planes. These give two disjoint embeddings $M_{\pm}=p_{\pm}(M)\hookrightarrow P_-$\,.\\  
\indent
As $\dim P_-=5$, $\dim\F=2$ and $\dim M_{\pm}=4$\,, at each null plane $p\in M_{\pm}$ we have
$d_{\pm}(p)=\dim(\F_p\cap T_pM_{\pm})\in\{1, 2\}$\,.\\ 
\indent
By Proposition \ref{prop:Grassmann}, $d_{\pm}=2$ on $M_{\pm}$ if and only if 
the two distributions $p_+$ and $p_-$ are integrable and totally geodesic. This is equivalent 
to the 
fact that the almost Hermitian structures $J_{\pm}$ defined by $J_{\pm}(X)=\pm\,{\rm i}X$ for 
$X\in p_{\pm}$ are integrable. Then, $\phi$ maps any leaf of $p_+$ and $p_-$ to a null 
geodesic on $N^2$\,; in particular, $\phi$ is horizontally conformal. Thus, a necessary 
condition for $\phi$ to be $(-)$twistorial is that $\phi$ be horizontally conformal. 
Then, on setting  
$P'=M_+\sqcup M_-$ and $\Phi=\dif\!\phi|_{P'}:P'\to Q$\,,  
\emph{$\phi$ is a\/ $(-)$twistorial map with respect to\/ $\Phi$ if and only if\/ $J_-(=-J_+)$ 
is integrable.}\\ 
\indent 
If $\phi$ is $(-)$twistorial, then the anti-self-dual 
surfaces on $M^4$ which are leaves of the distributions $p_+$ or $p_-$ are foliated by 
horizontal null geodesics (this follows from the fact that $\phi$ is horizontally conformal 
and, by definition, at each point $x\in M^4$, the space $p_{\pm}(x)$ intersects the horizontal  space of $\phi$ at $x$ along a null direction). 
If, further, $M^4$ is self-dual, then $Z(M,\phi)=S_+\sqcup S_-$ is the leaf space of the 
foliation induced by $\F$ on $M_+\sqcup M_-$, and the map $Z(M,\phi)\to Z(M)$ is simply 
the inclusion $S_+\sqcup S_-\hookrightarrow Z(M)$\,. Obviously $S_+$ and $S_-$ are 
transversal to the twistor lines. Furthermore, as the preimage of a null geodesic on $N^2$ 
through $\phi$ is a hypersurface on $M^4$ foliated by anti-self-dual surfaces, the two `surfaces' 
$S_+$ and $S_-$ in $Z(M)$ are foliated by curves with leaf spaces $C_+$ and $C_-$. 
Also, $Z(N,\phi)=Z(N)$ is the space of null geodesics of $N^2$ which, as a consequence of the 
horizontal conformality of $\phi$\,, can be canonically identified with $C_+\sqcup C_-$\,; then 
the twistorial representation of $\phi$ is simply the projection
$S_+\sqcup S_-\to Z(N)$. To retrieve $\phi$ from its twistorial representation, let $x\in M^4$;  
this corresponds to the twistor line $t_x\subseteq Z(M)$. Then, locally, $t_x$ meets $S_+$ and 
$S_-$ in two points which are projected by $S_+\sqcup S_-\to Z(N)$ to two points: one in 
$C_+$ and the other in $C_-$. These correspond to two null geodesics on $N^2$ whose intersection is precisely $\phi(x)$.\\ 
\indent
In a similar way to Example \ref{exm:3to2}, if $\phi:M^4\to N^2$ is a horizontally conformal 
submersion for which $J_{-}$($=-J_+$) is nowhere integrable then $d_{\pm}=1$ on $M_{\pm}$ and, 
if $M^4$ is self-dual then, locally, each anti-self-dual surface contains precisely one horizontal 
null geodesic and so we have two \emph{local} sections $Z_{\pm}\hookrightarrow\n(M)$ of the 
canonical projection 
(see \cite{Hit-complexmfds}) $\n(M)\to Z(M)${\,, where $\n(M)$ is the space of null geodesics 
of $M^4$ \cite{LeB-nullgeod}\,}. Note that if $\phi$ is $(-)$twistorial then the 
fibres of $\n(M)\to Z(M)$ induce on $Z_+\sqcup Z_-$ a one-dimensional foliation whose leaf 
space is $S_+\sqcup S_-$\,.\\ 
\indent 
A horizontally conformal submersion $\phi:M^4\to N^2$ with nowhere degenerate fibres is 
$(-)$twistorial if and only if its fibres are $(-)$twistorial in the sense of Example \ref{exm:2to4}\,. 
In particular, \emph{a horizontally conformal submersion $\phi:M^4\to N^2$ with nowhere degenerate fibres 
is both $(+)$twistorial and $(-)$twistorial if and only if it has totally umbilical fibres} 
(cf.\ Example \ref{exm:2to4}).  
\end{exm} 

\indent
The following theorem is a reformulation of a result of \cite{Woo-4d} 
(cf.\ Proposition \ref{prop:3to2harmorph} and Theorem \ref{thm:4to3harmorph}, below).  

\begin{thm} \label{thm:4to2harmorph} 
Let $M^4$ be an orientable four-dimensional Einstein manifold and let $\phi:M^4\to N^2$ be 
the complexification of a submersive harmonic morphism.\\ 
\indent
Then, with a suitable choice of orientations, $\phi$ is $(-)$twistorial. \qed
\end{thm}

\begin{rem} 
1) Let $M$ be a manifold endowed with an almost twistorial structure $\t=(P,M,\p,\F)$, $(\dim\F=k)$, 
such that the map $P\to G_k(TM)$ given by $p\mapsto\dif\!\p(\F_p)$ is an embedding of $P$ in the Grassmann 
bundle of $k$-dimensional planes on $M$.\\  
\indent
Then $\t$ is integrable if, locally, there are sufficiently many twistorial 
maps from $(M,\t)$; for example, if for each $p\in P$ there exists a (local) twistorial map 
$\phi$ from $(M,\t)$ such that ${\rm ker}\dif\!\phi_{\p(p)}=p$ then $\t$ is integrable. In the above examples 
this condition is also necessary; moreover other similar necessary and sufficient conditions can easily 
be formulated. In Theorem \ref{thm:4to3} below, we shall see that the existence of a \emph{single} suitable twistorial map may be sufficient for a twistorial structure to be integrable.\\   
\indent 
2) Twistorial maps as in Example \ref{exm:3to2}, Example \ref{exm:4to2} and as in the next 
section appear, in a more or less explicit way, in \cite{LeB-jdg}. 
\end{rem}

\section{Twistorial maps with one-dimensional fibres from four-dimensional 
Riemannian manifolds} \label{section:4to3}

\indent
In this section we continue to work in the complex-analytic category. 
The almost twistorial structures which may appear in this section will be those of 
Examples \ref{exm:3t} and \ref{exm:4t}.
As usual, the results can be applied to the real-analytic category by complexification (Remark \ref{rem:tm}).
\\ 
\indent
Let $\phi:M^4\to N^3$ be a surjective submersive map with nowhere degenerate fibres from  an oriented $4$-dimensional Riemannian manifold to a $3$-dimensional Riemannian manifold . 
Suppose, for the moment, that $M^4$ is self-dual and is endowed with the twistorial structure 
of Example \ref{exm:4t}(ii)\,, with $P=P_-$ the bundle of anti-self-dual planes of $M^4$, and that 
$N^3$ is endowed with the twistorial structure of Example \ref{exm:3t}\,, which corresponds to an 
Einstein--Weyl connection $D$ on $N^3$, with $Q$ the bundle of degenerate planes on $N^3$.  
Then, it is clear that $\phi$ is $(-)$twistorial for a suitable choice of 
$P' \subseteq P$ if and only if it maps some, if not all, of the anti-self-dual surfaces on $M^4$ to totally geodesic degenerate  
surfaces on $(N^3,D)$ (cf.\ \cite{Hit-complexmfds}, \cite{JonTod}). But, unless we introduce 
some extra structure, there is no reason to ignore any of the anti-self-dual surfaces on $M^4$. 
Moreover, if $\phi$ maps anti-self-dual surfaces on $M^4$ to totally geodesic 
degenerate surfaces on $(N^3,D)$, then, in particular, $\dif\!\phi$ maps anti-self-dual planes on $M^4$ to degenerate planes on $N^3$ which is equivalent to the condition that
$\phi$ be horizontally conformal (indeed, the differential $\dif\!\phi_x$\,, 
($x\in M$), maps the orthogonal complement, in $\H_x$\,, of a null direction 
$l\subseteq \H_x$ onto the orthogonal complement of a null direction 
$k\subseteq T_{\phi(x)}N$ if and only if the adjoint of $\dif\!\phi_x$ maps 
$k$ onto $l$; thus, as the horizontal projection of any anti-self-dual plane is a 
degenerate horizontal plane and any degenerate horizontal plane is obtained in this 
way, $\dif\!\phi$ maps anti-self-dual planes on $M^4$ to degenerate planes on $N^3$ 
if and only if the adjoint of $\dif\!\phi$ is conformal at each point).\\
\indent
Therefore, given an oriented Riemannian manifold $M^4$ and a 
horizontally conformal submersion $\phi:M^4\to N^3$ with nowhere degenerate fibres onto 
a Riemannian manifold, we shall look for necessary and sufficient conditions under which 
there exists a Weyl connection $D$ on $N^3$ with respect to which $\phi:M^4\to(N^3,D)$ is a 
$(-)$twistorial map with $P'=P_-$ the bundle of anti-self-dual planes on $M^4$.\\ 
                                                                                                
\indent
To do this, we first place the discussion in a slightly more general context. 
Let $M^4$ be an oriented 
Riemannian manifold and let $\H$ be a nowhere degenerate three-dimensional distribution. Denote, as 
usual, $\V=\H^{\perp}$ and assume that $(\V,g|_{\V})$ and $(\H,g|_{\H})$ are oriented so that 
the isomorphism $(TM,g)=(\V,g|_{\V})\oplus(\H,g|_{\H})$ is orientation preserving. Let $I^{\H}$ 
be the integrability 2-form of $\H$ defined by $I^{\H}(X,Y)=-g(U,[X,Y])$ 
for all local sections $X,Y$ of $\H$ where $U$ is the positive unit section of $\V$. 
Also, let $\Bh$ be the second fundamental form of $\H$ defined by 
$\Bh(X,Y)=\frac12\V(\nabla_XY+\nabla_YX)$ for any local horizontal vector fields 
$X$ and $Y$, where $\nabla$ is the Levi-Civita connection of $(M^4,g)$\,. Then we define 
a section $\UBh$ of $\H^*\otimes\H^*$ by 
$\UBh(X,Y)=g(U,\Bh(X,Y))$ for any horizontal $X$ and $Y$; we shall denote by $\UB0h$ the trace-free part of $\UBh$\,.  Let $*_{\H,g}$ be the Hodge star-operator on 
$(\H,g|_{\H})$\,. 

\indent 
Recall (see \cite{Gau}\,) that if $D$ is a conformal connection on a conformal manifold $(M,c)$ 
and $g$ is a local representative of $c$ on some open set $U$, then the \emph{Lee form} of $D$ 
\emph{with respect to $g$} is the one-form $\a\in\G(T^*U)$ characterised by $Dg=-2\a\otimes g$\,. 
The Lee form of a conformal \emph{partial} connection can be defined in a similar way. 
Also, if $D$ is a partial connection 
on $\H$, over $\H$\,, then its \emph{torsion} $T$\,, with respect to $\V$, is defined by 
$T(X,Y)=D_XY-D_YX-\H[X,Y]$ for any local sections $X$ and $Y$ of $\H$\,. We now introduce a
conformal partial connection which encodes the condition of twistoriality.
 
\begin{defn}[cf.\ \cite{Cal-sds}] \label{defn:Dpm}
The \emph{conformal partial connections $D_{\pm}$ induced by $g$ on $\H$} 
are the unique conformal partial connections on $(\H,g|_{\H})$\,, over $\H$\,, whose Lee forms 
$\a_{\pm}$ and torsion tensors $T_{\pm}$\,, with respect to $\V$, are given by 
\begin{equation} \label{e:t4to3}
\begin{split}
\a_{\pm}&=\trace(\Bv)^{\flat}\pm*_{\H,g}I^{\H}\;,\\
*_{\H,g}T_{\pm}&=\mp2\UB0h
\end{split}
\end{equation}
where we have identified $\H\otimes\H^*$ and $\H^*\otimes\H^*$ by using $g|_{\H}$\,. 
\end{defn}

\indent 
If $\H$ is totally umbilical then $D_-$ is the horizontal part of a Weyl connection of 
\cite[\S4]{Cal-sds}\,.\\

\indent 
Let $\p^{\H}:P^{\H}\to M$ be the bundle of degenerate planes of $\H$\,. As in Example \ref{exm:3t}\,, $D_-$  
induces a twistor distribution $\F_-^{\!\H}$ defined by setting 
$\F_-^{\!\H}(p)\subseteq T_pP^{\H}$ equal to the horizontal
lift at $p\in P^{\H}$, with respect to $D_-$\,,  of $p\subseteq\H_{\p^{\H}(p)}$. 
Assume, for simplicity, that $\V={\rm ker}\dif\!\phi$ where $\phi:M^4\to N^3$ is a 
surjective submersion.

\begin{prop} \label{prop:Dpm} 
{\rm (a)} Let $\g$ be a null geodesic on $(M^4,g)$\,.\\ 
\indent 
{\rm (1)} If $\g$ is horizontal then it is a geodesic of\/ $D_-$\,.\\ 
\indent 
{\rm (2)} Suppose that $\g$ is not horizontal at $x_0=\g(0)$\,;  
let $p_{0-}\subseteq T_{x_0}M$ be the anti-self-dual plane containing $(\dif\!\g/\dif\!z)(0)$ 
and let $p_0=\H(p_{0-})$\,. 
Denote by\/ $c$ the (local) horizontal curve such that $c(0)=x_0$ and 
$\phi\circ c=\phi\circ\g$\,.\\ 
\indent 
If $p$ is the field of horizontal degenerate planes along $c$ such that $p(0)=p_0$ then 
$(\dif\!p/\dif\!z)(0)\in\F_-^{\!\H}$\,.\\ 
\indent 
{\rm (b)} Moreover, if $D$ is a conformal partial connection on $(\H,g|_{\H})$ which has this property 
for all nonhorizontal null geodesics on $(M^4,g)$\,, and whose torsion $T^D$\,, with respect to $\V$, 
is such that $*T^D$ is self-adjoint and trace-free, then $D=D_-$\,. In particular, $D_-$ is conformally invariant. 
\end{prop} 
\begin{proof} 
(1) Suppose that $\g$ is horizontal and let $Y$ be its velocity vector field. Then 
$\UB0h(Y,Y)=0$ and hence 
also $(*_{\H,g}T_{-})(Y)=0$\,; equivalently, $T_-$ restricted to $Y^{\perp}$ is zero, and thus 
$g(Y,T_-(X,Y))=0$ for any $X\in Y^{\perp}$.\\ 
\indent 
Because $\g$ is a null curve and $D\!_-$ is conformal, we have that 
$D\!_{-}{}_{\;Y}Y\in Y^{\perp}$ 
at each point of $\g$\,. 
Thus, to prove that $\g$ is a geodesic of $D\!_-$\,, it is sufficient to prove 
that $g(D\!_{-}{}_{\;Y}Y,X)=0$ for any $X\in Y^{\perp}$. The proof follows {from the 
fact that $D\!_{-}{}_{\;Y}Y=\nabla_YY+2\a_-(Y)Y$ where $\nabla$ is the Levi-Civita 
connection of $g$\,}.\\ 
\indent 
(2) If $\g$ is not horizontal then we can write $\dif\!\g/\dif\!z=X+{\rm i}\,U$ 
where $X$ is horizontal, $U$ is vertical and $g(X,X)=g(U,U)=a^2\neq0$\,. Extend $X$ and $U$ to local 
sections of $\H$ and $\V$, respectively. We can assume that $X$ is a basic vector field and $a^{-1}U$ 
is positive on $(\V,g|_{\V})$. Then, along $\g$\,, we have 
\begin{equation} \label{e:trivi}  
\nabla_{X+{\rm i}\,U}(X+{\rm i}\,U)=\nabla_XX-\nabla_UU+{\rm i}\,(\nabla_XU+\nabla_UX)=0  
\end{equation} 
where $\nabla$ is the Levi-Civita connection of $(M^4,g)$\,.\\ 
\indent 
We see from \eqref{e:trivi} that, for any horizontal vector field $Y$, we have  
\begin{equation*} 
g(\nabla_XX,Y)-g(\nabla_UU,Y)+{\rm i}\,g(\nabla_XU,Y)+{\rm i}\,g(\nabla_UX,Y)=0   
\end{equation*} 
along $\g$\,.\\ 
\indent 
Now, assume that $Y$ is null (and horizontal) and $g(X,Y)=0$\,. Then, by using the fact that 
$X$ is basic, the last relation implies that 
\begin{equation} \label{e:trivia} 
g(\nabla_XX,Y)=a^2\bigl\{g(\trace(\Bv),Y)-{\rm i}\,I^{\H}(a^{-1}X,Y)+2{\rm i}\UBh(a^{-1}X,Y)\bigr\} 
\end{equation} 
along $\g$ where we have assumed that $\{a^{-1}U\}$ is a positive local frame of $(\V,g|_{\V})$\,.\\ 
\indent 
On the other hand, let $D$ be a conformal partial connection on $(\H,g|_{\H})$\,, let $T$ be its 
torsion with respect to $\V$, and let $\a$ be its Lee form with respect to $g|_{\H}$\,. Then 
we have  
\begin{equation} \label{e:aT} 
g(\nabla_XX,Y)=g(D_XX,Y)+\a(Y)g(X,X)+g(X,T(X,Y))\;. 
\end{equation} 
\indent 
Let $\{X_0,X_1,X_2,X_3\}$ be a positive local orthonormal frame on $(M^4,g)$ such 
that $X_0=a^{-1}U$ and $X_1=a^{-1}X$\,. Then  
$\bigl\{X_1,X_2,X_3\bigr\}$ is a positive orthonormal local frame of $(\H,g|_{\H})$. Take 
$Y=X_2+{\rm i}X_3$ 
and note that $p_{0-}$ is spanned by $\{X_{x_0}+{\rm i}\,U_{x_0},Y_{x_0}\}$ and $p_0$ by 
$\{X_{x_0},Y_{x_0}\}$\,. Also, $c$ is the integral curve of $X$ through $x_0$ so that $p_z$ is spanned 
by $\{X_{c(z)},Y_{c(z)}\}$ at each $z$\,.\\ 
\indent 
Then $(\dif\!p/\dif\!z)(0)$ is horizontal with respect to the connection induced by $D$ 
on $\p^{\H}:P^{\H}\to M$ if and only if $(D_XX)_{x_0}\in p_0$\,; equivalently  
$g\bigl((D_XX)_{x_0},Y_{x_0}\bigr)=0$\,. 
By \eqref{e:trivia} and \eqref{e:aT}\,, the last condition is equivalent to 
\begin{equation} \label{e:almostalmostthere} 
\a(Y)-g(\trace(\Bv),Y)+{\rm i}\,I^{\H}(X_1,Y)= 
-g(X_1,T(X_1,Y))+2{\rm i}\UBh(X_1,Y)  
\end{equation} 
at $x_0$\,. Condition \eqref{e:almostalmostthere} is equivalent to 
\begin{equation} \label{e:almostthere} 
\bigl(\a-\trace(\Bv)^{\flat}+*_{\H,g}I^{\H}\bigr)(Y)= 
-g(X_1,T(X_1,Y))+2{\rm i}\UBh(X_1,Y)
\end{equation} 
at $x_0$\,. 
It is easy to see that, if in \eqref{e:almostthere} we replace $\a$ and $T$ by $\a_-$ and $T_-$\,, 
respectively, then both sides are zero. Hence, $(\dif\!p/\dif\!z)(0)\in\F_-^{\!\H}$\,.\\ 
\indent 
Suppose now that $D$ satisfies condition (2) of the proposition. Then \eqref{e:almostthere} must hold. Moreover, the relation obtained from \eqref{e:almostthere} by replacing $Y$ with
$\widetilde{Y}=X_2-{\rm i}X_3$ and $X_1$ by $-X_1$ (so that $\{X_0,-X_1,X_2,-X_3\}$ is positive) must also hold, viz. 
\begin{equation} \label{e:nearlythere} 
\bigl(\a-\trace(\Bv)^{\flat}+*_{\H,g}I^{\H}\bigr)(\widetilde{Y})= 
-g(X_1,T(X_1,\widetilde{Y}))-2{\rm i}\UBh(X_1,\widetilde{Y})
\end{equation} 
at $x_0$\,. By taking the sum of \eqref{e:almostthere} and \eqref{e:nearlythere} we obtain  
\begin{equation} \label{e:closelythere} 
\bigl(\a-\trace(\Bv)^{\flat}+*_{\H,g}I^{\H}\bigr)(X_2)= 
-g(X_1,(*T)(X_3))-2\UBh(X_1,X_3) 
\end{equation} 
at $x_0$\,.\\ 
\indent 
{}From the fact that the right-hand-side of \eqref{e:closelythere} does not depend on $X_2$ it follows 
that $*T+2\UBh$ is proportional to $g|_{\H}$\,. But $*T$ is self-adjoint and trace-free, hence 
$*T+2\UB0h=0$\,. Then, by \eqref{e:closelythere}, we also have  
$\a=\trace(\Bv)^{\flat}-*_{\H,g}I^{\H}$ and hence $D=D_-$\,.
\end{proof} 

\indent
Since, up to Riemannian coverings, any horizontally conformal 
submersion admits global dilations, from now on, for simplicity, we shall consider only 
horizontally conformal submersions which admit global dilations.\\ 
\indent
The following theorem is a reformulation of results of \cite{Cal-sds} (see Definitions \ref{defn:(-)twist}
and \ref{defn:Dpm} for `$(-)$twistorial map' and $D_-$\,, respectively).
   
\begin{thm} \label{thm:4to3} 
Let $(M^4,g)$ and $(N^3,h)$ be orientable Riemannian manifolds and let $\phi:(M^4,g)\to(N^3,h)$ be a 
surjective horizontally conformal submersion with connected nowhere degenerate fibres.\\ 
\indent 
Suppose that orientations on $(M^4,g)$\,, $(\V,g|_{\V})$ and $(\H,g|_{\H})$ are chosen 
such that the isomorphism $(TM,g)=(\V,g|_{\V})\oplus(\H,g|_{\H})$ is orientation preserving.\\ 
\indent   
Then, the following assertions are equivalent:\\ 
\indent
{\rm (i)} there exists a Weyl connection $D$ on $N^3$ with respect to which $\phi$ is $(-)$twistorial (with $P'=P_-$ the bundle of anti-self-dual planes on $(M^4,g)$\,);\\ 
\indent
{\rm (ii)} the one-form
\begin{equation} \label{e:tl4to3}
\trace(\Bv)^{\flat}-\l^{-1}\dif^{\H}\!\l-*_{\H,g}I^{\H}
\end{equation}
is basic where $\l$ is the dilation of\/ $\phi$, $\dif^{\H}=\H\circ\dif$ and $*_{\H,g}$ is the Hodge
star-operator on $(\H,g|_{\H})$\,;\\ 
\indent
{\rm (iii)} the two-form $\dif\bigl(\trace(\Bv)^{\flat}+\tfrac13\,\trace(\Bh)^{\flat}\bigr)$ is 
self-dual.\\  
\indent
Moreover, if any of the assertions {\rm (i)}, {\rm (ii)} or {\rm (iii)} holds then:\\ 
\indent
{\rm (a)} $\phi^*(D)^{\H}=D_-$\,; in particular, $D$ is the \emph{unique} Weyl connection that satisfies 
{\rm (i)}\,;\\  
\indent
{\rm (b)} $M^4$ is self-dual if and only if $D$ is Einstein--Weyl.\\  
\indent
If $M^4$ is self-dual then, the following assertion can be added to the assertions {\rm (i)} to {\rm (iii)}:\\ 
\indent
{\rm (iv)} There exists a (unique) Einstein--Weyl connection $D$ on $N^3$ such that, for any local 
map $\psi:N^3\to P^2$ which is twistorial in the sense of Example \ref{exm:3to2}\,, 
the map $\psi\circ\phi:M^4\to P^2$ is $(-)$twistorial 
(in the sense of Example \ref{exm:4to2}\,).  
\end{thm}  
\begin{proof} 
Let $\p_-:P_-\to M$ be the bundle of anti-self-dual planes on $M^4$ and let $\Phi=\dif\!\phi|_{P_-}$ 
be the map induced by $\phi$ from $P_-$ to the bundle of degenerate planes on $N^3$. 
Assertion (i) is equivalent to the existence of a Weyl connection $D$ on $N^3$ such that for any $x\in M^4$ and any null geodesic $\g$ on $M^4$ with $\g(0)=x$\,, if we denote 
by $p$ the field of anti-self-dual planes along $\g$  
such that $(\dif\!\g/\dif\!z)(z)\in p(z)$ for all $z$\,, then 
$\dif\!\Phi\bigl(\dif\!p/\dif\!z)(0)\bigr)\in\F$\,, where $\F$ is the twistor distribution induced by $D$ on the bundle of degenerate planes on $N^3$ (Example \ref{exm:3t}). 
Thus, by Proposition \ref{prop:Dpm}\,, assertion (i) is equivalent 
to the existence of a Weyl connection $D$ on $N^3$ such that $\phi^*(D)^{\H}=D_-$\,.  
Now, applying the conformal change $\phi^*(h)|_{\H}=\l^2\,g|_{\H}$ to \eqref{e:t4to3}
introduces the term $-\l^{-1}\!\dif^{\H}\!\l$ so that the Lee form of $D_-$ with respect
to $\phi^*(h)|_{\H}$ is equal to the one-form \eqref{e:tl4to3}\,.
Hence we have the equivalence of (i) and (ii)\,. Furthermore, if (i) holds then assertion
(a) also holds.\\ 
\indent 
The equivalence (ii)$\iff$(iii) can be found in the proof of 
\cite[Proposition 4.4]{Cal-sds}\,.\\ 
\indent
The fact that, if $M^4$ is self-dual, then (i)$\iff$(iv) is obvious. Also, if $M^4$ is self-dual 
and $\phi$ is $(-)$twistorial, then from Theorem \ref{thm:3t} it follows that $D$ is Einstein--Weyl.  
Conversely, if $\phi$ is $(-)$twistorial and $D$ is Einstein--Weyl, then by composing $\phi$ with local twistorial maps $(N^3,D)\to P^2$ we obtain sufficiently many local 
twistorial maps $M^4\to P^2$ to show that any anti-self-dual plane on $M^4$ is tangent 
to some anti-self-dual surface in $M^4$. Thus, by Theorem \ref{thm:PenAtHiSi-4t}\,, 
$M^4$ is self-dual. 
\end{proof} 

\begin{rem} 
1) D.M.J.~Calderbank \cite{Cal-sds} calls a horizontally conformal submersion, between \emph{real-}analytic
Riemannian manifolds, satisfying condition (iii) of Theorem \ref{thm:4to3} \emph{self-dual}. We have thus
interpreted self-dual submersions as those which are $(-)$twistorial.\\
\indent 
2) By using the null-tetrad formalism (see \cite{MasWoo}\,), the equivalence (i)$\iff$(ii) of 
Theorem \ref{thm:4to3} can also be obtained after a straightforward but tedious computation.  
\end{rem} 

\indent 
In the following, let $\phi:(M^4,c_M)\to(N^3,c_N)$ be a surjective horizontally conformal submersion 
with nowhere degenerate fibres from a four-dimensional 
oriented conformal manifold to a three-dimensional conformal manifold. Let $L$ be the line bundle over 
$N^3$ associated with the bundle of conformal frames on $N^3$ through the morphism of Lie groups 
$\rho:CO(3, \mathbb{K})\to\mathbb{K}\setminus\{0\}$\,, 
($\mathbb{K}=\R$ or $\C$)\,, characterised by $a\in CO(3,\mathbb{K})$ 
if and only if, $a^t\,a=\rho(a)^2I$
(in the smooth category, $L$ can be defined as the line bundle associated to the frame bundle of $N^3$ through 
the morphism of Lie groups $GL(3,\R)\to(0,\infty)$\,, $a\mapsto|\det(a)|^{1/3}$\,, see \cite{Cal-sds}\,). 
Note that $L$ can be defined on any odd-dimensional conformal manifold (in the smooth category, 
on \emph{any} (conformal) manifold) $(N,c_N)$ and that its local sections 
correspond to oriented local representatives of $c_N$\,. 
Let $E$ be some line bundle over $N$ and let $\widetilde{E}=\phi^*(E)$\,.\\ 
\indent 
To any pair $(s,\nabla)$ where $s$ is a section of $L^*$ over some open subset $U$ of $N$ 
and $\nabla$ is a connection on $E|_U$ there can be associated a connection $\widetilde{\nabla}$ of  
$\widetilde{E}|_{\phi^{-1}(U)}$ as follows: assume, initially, that $s$ is nowhere zero, let $h$ 
be the oriented local representative of $c_N$ over $U$  
corresponding to it and let $g$ be the oriented local representative of $c_M$ over 
$\phi^{-1}(U)$ 
such that $\phi:(\phi^{-1}(U),g)\to(U,h)$ is a Riemannian submersion.  Denote, as 
usual, $\V={\rm ker}\dif\!\phi$\,, $\H=\V^{\perp}$, and let $\omega\in\G(\phi^{-1}(U),\V^*)$ be the 
induced orientation on $\V|_{\phi^{-1}(U)}$ such that 
the isomorphism $\bigl(T(\phi^{-1}(U)),g\bigr)=\bigl(\V|_{\phi^{-1}(U)},g|_{\V}\bigr)\oplus\bigl(\H|_{\phi^{-1}(U)},\phi^*(h)\bigr)$ 
is orientation preserving. 
If $A$ is the local connection form of $\nabla$ with respect to some nowhere zero 
local section $\s\in\G(U,E)$\,,  
then we set the local connection form of $\widetilde{\nabla}$ corresponding to $\phi^*(\s)$ equal to 
$$\widetilde{A}=-\omega+\phi^*(A)\;.$$ 
We shall call $\widetilde{\nabla}$ the \emph{pull-back, by $\phi$\,, of the pair $(s,\nabla)$\,.} Note that, if $u$ is a nowhere zero function locally defined on $N$, 
the pull-back by $\phi$ of $(us,\nabla)$ is given by $-u\,\omega+\phi^*(A)$\,. Hence, we can extend the 
definition of the pull-back by $\phi$ of any pair $(s,\nabla)$ where the local section 
$s$ of $L^*$ may have zeros.\\ 
\indent 
The monopole equation of \cite{JonTod} can be written as follows:
\begin{equation} \label{e:JonTod}
(\dif -\a)u=*F
\end{equation}
where $u$ is a function on a three-dimensional oriented Riemannian manifold $(N^3,h)$,
$\a$ is a one-form on $N^3$, $F$ is a two-form on $N^3$ and $*$ is the Hodge star-operator
of $(N^3,h)$\,. By interpreting $\a$ as the Lee form with respect to $h$ of a Weyl connection
on $(N^3,[h])$, and $F$ as the curvature form of some connection on some line bundle over $(N^3,h)$,
the equation \eqref{e:JonTod} can be written as follows:

\begin{defn}[\cite{JonTod,CalPed}] \label{defn:JonTodCalPed}
Let $D$ be a Weyl connection on $(N^3,c_N)$\,. A \emph{monopole on $(N^3,c,D)$} is a pair 
$(s,\nabla)$ where $s$ is a section of $L^*$ and $\nabla$ is a connection on a line bundle $E$ over $N^3$ such that 
\begin{equation} \label{e:JonTodCalPed} 
*_NDs=F
\end{equation}  
where $F$ is the curvature form of $\nabla$\,. A monopole is called \emph{nontrivial} if $s\neq0$\,. 
\end{defn} 

\begin{rem}  
Let $(N^3,c)$ be a three-dimensional conformal manifold endowed with a Weyl connection $D$\,.\\ 
\indent 
1) It is well known (see \cite{CalPed}\,) that if $(N^3,c,D)$ is Einstein--Weyl then, 
at least locally, there exist nontrivial monopoles on $(N^3,c,D)$\,.\\  
\indent 
Conversely, let $(N^3,c)$ be a three-dimensional conformal manifold and let $E$ be a 
line bundle over $N$ endowed with a connection $\nabla$. Also, let $s$ be a nowhere 
zero (local) section of $L^*$ and let $h$ be the oriented (local) representative of 
$c$ corresponding to $s$\,. Define a one-form $\a$ such that $-*_h\a=\phi^*(F)$ where 
$F$ is the curvature form of $\nabla$\,. Then $(s,\nabla)$ is a monopole on $(N^3,c,D)$ 
where $D$ is the Weyl connection on $(N^3,c)$ whose Lee form with respect to $h$ is $\a$\,.\\  
\indent 
2) Let $h$ be a local representative of $c$ over some open subset $U$ of $N$ and let 
$\a$ be the Lee form of $D$ with respect to $h$\,. Then, it is well known (and easy 
to prove) that, on small enough open subsets $U$\,, monopoles on $(N^3,c,D)$ correspond 
to solutions of the equation $\D u=\dif^*(u\a)$ where 
$\D$ is the Laplacian and $\dif^*$ is the codifferential on $(U,h)$\,. 
Hence, if $(N^3,c,D)$ is (the complexification) of a real-analytic conformal manifold 
endowed with a real-analytic Weyl connection then, at least locally, there exist nontrivial monopoles on 
$(N^3,c,D)$\,. 
\end{rem} 

\indent 
Assertion (i) of the following proposition is essentially in \cite{CalPed,Cal-sds}\,.

\begin{cor} \label{cor:monomor} 
Let $\phi:(M^4,c_M)\to(N^3,c_N)$ be a surjective horizontally conformal submersion with 
nowhere degenerate 
fibres and let $D$ be a Weyl connection on $(N^3,c_N)$\,.\\ 
\indent 
{\rm (i)} If $\phi:(M^4,c_M)\to(N^3,c_N,D)$ is a $(-)$twistorial map, then it pulls back 
monopoles on $(N^3,c_N,D)$ to self-dual connections on $(M^4,c_M)$\,.\\  
\indent 
{\rm (ii)} Conversely, suppose that there exists a nontrivial monopole on $(N^3,c_N,D)$ which is pulled back 
by $\phi$ to a self-dual connection on $(M^4,c_M)$\,. Then, $\phi:(M^4,c_M)\to(N^3,c_N,D)$ is $(-)$twistorial.  
\end{cor} 
\begin{proof} 
Let $(s,\nabla)$ be a nontrivial monopole. We may assume that $s$ is nowhere zero so that it corresponds 
to an oriented local representative $h$ of $c_N$\,. Let $\a$ be the Lee form of $D$ with respect to 
$h$\,. Then \eqref{e:JonTodCalPed} is equivalent to 
\begin{equation} \label{e:JTCP} 
-*_h\a=F\;. 
\end{equation} 
\indent 
If $A$ is a local connection of $\nabla$\,, then the corresponding local connection form of $\widetilde{\nabla}$ 
is $\widetilde{A}=-\omega+\phi^*(A)$ where $\omega\in\G(\V^*)$ is as above. Now, by using the fact that 
\begin{equation} \label{e:difomega}
\dif\!\omega=-\trace(\Bv)^{\flat}\wedge\omega+I^{\H} 
\end{equation} 
it is easy to prove that
$\dif\!\widetilde{A}$ is self-dual if and only if 
\begin{equation} \label{e:JTCP-} 
-*_h\a_-=\phi^*(\dif\!A)  
\end{equation}   
where $\a_-=\trace(\Bv)^{\flat}-*_{\H,g}I^{\H}$ is the Lee form of $D_-$ with respect to $h$\,.\\
\indent 
{}From \eqref{e:JTCP} and \eqref{e:JTCP-} it follows that $\a_-=\a$ is basic. The result follows from
Theorem \ref{thm:4to3}\,.
\end{proof}

\indent 
Let $M^m$ be a manifold of dimension $m\geq3$ endowed with a three-dimensional distribution $\H$\,. 
Suppose that $\H$ is endowed with a conformal structure $c$ and a conformal partial connection $D$\,. 
In a similar fashion to the above, we consider a line bundle $L$ over $M$ whose nowhere zero local 
sections correspond to oriented local representatives of $c$\,.\\ 
\indent  
Suppose that there exists a foliation $\V$ which is complementary to $\H$ and assume, for simplicity,  
that there exists a submersion $\phi:M^m\to N^3$ whose fibres are leaves of $\V$.\\ 
\indent 
Then Definition \ref{defn:JonTodCalPed} can be easily generalised by defining a 
\emph{monopole on $(\H,c,D)$ (with respect to $\phi$)}  
to be a pair $(s,\nabla)$ where $s$ is a section of $L^*$ and $\nabla$ a connection on a line bundle $E$ over $N^3$ such that 
$*_{\H}Ds=\phi^*(F)$ where $F$ is the curvature form of $\nabla$\,. (One could also consider line 
bundles on $M$ endowed with partial connections over $\H$ but then the resulting equation would depend  
on a local section of $E$\,.)\\  
\indent 
Furthermore, if $(M^4,c)$ is a four-dimensional oriented conformal manifold then, a construction similar 
to above associates to any pair $(s,\nabla)$ 
a connection $\widetilde{\nabla}$ on $\phi^*(E)$\,. We then have the following simple fact (which may be 
known) whose proof is similar to the proof of Corollary \ref{cor:monomor}\,.  
 
\begin{prop} \label{prop:mono-sd} 
Let $\phi:(M^4,c)\to N^3$ be a surjective submersion with nowhere degenerate fibres from an oriented 
four-dimensional conformal manifold. Suppose that $D$ is a conformal partial connection on $(\H,c|_{\H})$ 
and let $E$ be a line bundle over $N^3$\,. Let $s$ be a section of\/ $L^*$ and $\nabla$ a connection on $E$\,.\\ 
\indent 
Then any two of the following assertions imply the third:\\ 
\indent 
{\rm (i)} $(s,\nabla)$ is a nontrivial monopole on $(\H,c|_{\H},D)$\,;\\ 
\indent 
{\rm (ii)} $\widetilde{\nabla}$ is self-dual;\\ 
\indent 
{\rm (iii)} $D=D_-$\,. \qed
\end{prop}  

\indent 
Next, we show how harmonic morphisms fit into the above discussion. 

\begin{prop} \label{prop:harmono} 
Let $\V$ be (the complexification) of a one-dimensional (nowhere degenerate) conformal foliation on a 
four-dimensional oriented conformal manifold $(M^4,c)$\,; let $\H=\V^{\perp}$.\\ 
\indent 
Then, the following assertions are equivalent:\\ 
\indent 
{\rm (i)}  there exist local representatives $g$ of\/ $c$ with respect to which $\V$ is locally defined 
by harmonic morphisms;\\ 
\indent 
{\rm (ii)} there exist nontrivial monopoles $(s,\nabla)$ on $(\H,c|_{\H},D_-)$ for which  
$\widetilde{\nabla}$ is flat. 
\end{prop}     
\begin{proof} 
\indent  
By Theorem \ref{thm:Bry}\,, assertion (i) is equivalent to the following:\\ 
\indent 
(i$'$) locally, there exist representatives $g$ of\/ $c$ with respect to which 
$\V$ is geodesic.\\ 
\indent 
Also, from \eqref{e:difomega}\,, we easily obtain that a pair $(s,\nabla)$ is as 
in assertion (ii) if and only if the following hold:\\  
\indent 
(a) $\V$ has geodesic fibres with respect to the local representative $g$ of\/ $c$ whose 
horizontal component corresponds to $s$\,;\\ 
\indent 
(b) if $\theta$ is the local volume form induced by $g$ on $\V$ then $\dif\theta$ is 
equal to the pull-back by $\phi$ of the curvature form of\/ $\nabla$\,.\\ 
\indent 
We have thus shown that (ii)$\Longrightarrow$(i)\,.\\ 
\indent 
Conversely, suppose that (i) holds and let $s$ be a nowhere zero local section 
of\/ $L^*$ such that, with respect to the corresponding oriented local representative $g$  
of\/ $c$\,, the foliation $\V$ is geodesic. Then $I^{\H}=\dif\theta$ where 
$\theta$ is the local volume form induced by $g$ on $\V$\,. Moreover, because $\V$ 
is geodesic with respect to $g$\,, we have, locally, $\dif\theta=\phi^*(\dif\!A)$ for some 
one-form $A$ on $N^3$. Furthermore, by Definition \ref{defn:Dpm}\,, we have 
$-*_{\H,g}\a_-=I^{\H}$\,. Hence $-*_{\H,g}\a_-=\phi^*(\dif\!A)$\,; equivalently, the pair 
formed by $s$ and the connection determined by $A$\,, on the trivial line bundle 
over $N^3$, is a monopole on $(\H,c|_{\H},D_-)$\,. The proposition is proved. 
\end{proof} 

\indent 
The following result gives another characterisation for twistorial harmonic morphisms between orientable 
Riemannian manifolds of dimension four and three, respectively.  
  
\begin{cor} \label{cor:4to3harmorph} 
Let $(M^4,g)$ and $(N^3,h)$ be orientable Riemannian manifolds. Let $\phi:(M^4,g)\to(N^3,h)$ be a 
surjective submersive harmonic morphism with connected nowhere degenerate fibres. 
Suppose that orientations on $(M^4,g)$\,, $(N^3,h)$\,, $(\V,g|_{\V})$ and $(\H,g|_{\H})$ are chosen 
such that the isomorphisms $(TM,g)=(\V,g|_{\V})\oplus(\H|_{\H})$ and $(\H,g|_{\H})=(\phi^*(TN),\phi^*(h)\,)$ 
are orientation preserving.\\ 
\indent  
Then, the following assertions are equivalent:\\ 
\indent
{\rm (i)} There exists a Weyl connection $D$ on $N^3$ with respect to which $\phi$ 
is $(-)$twistorial.\\ 
\indent 
{\rm (ii)} There exists a basic one-form $\a$ on $M^4$ such that 
\begin{equation} \label{e:m} 
(\dif\!^{\H}-\a)(\l^{-2})=*_\H\O 
\end{equation} 
where $\l$ is a dilation of\/ $\phi$\,, $*_\H$ is the Hodge star operator on\/ $\bigl(\H,\phi^*(h)\bigr)$  
and\/ $\O$ is the curvature form of the horizontal distribution (i.e., in the notation of 
Theorem \ref{thm:Bry}\,, $\O=\dif\theta$).\\ 
\indent
Moreover, if {\rm (i)} or {\rm (ii)} holds then the Weyl connection $D$ of assertion {\rm (i)} is unique 
and the one-form $\a$ of assertion {\rm (ii)} is the pull-back of the Lee form of\/ $D$ with respect to 
$h$\,. 
\end{cor} 
\begin{proof} 
Comparing the fundamental equation \eqref{e:fundamentaleqn} with \eqref{e:tl4to3} and by using the fact that
$I^{\H}=\l\,\Omega$\,, it follows that a
one-form $\a$ on $M^4$ satisfies \eqref{e:m} if and only if it is the Lee form of\/ $D_-$ with respect to
$\phi^*(h)$\,. The proof follows from Theorem \ref{thm:4to3}\,.
\end{proof} 

\indent 
{}From \cite[Corollary 1.9]{Pan-4to3} we obtain the following result on maps from Einstein manifolds 
(cf.\ Proposition \ref{prop:3to2harmorph} and Theorem \ref{thm:4to2harmorph}). 

\begin{thm} \label{thm:4to3harmorph} 
Let $M^4$ be an orientable four-dimensional Einstein manifold and let $\phi:M^4\to N^3$ be 
the complexification of a submersive harmonic morphism.\\ 
\indent
Then, there exists a Weyl connection on $N^3$ with respect to which, with a suitable 
choice of orientations, $\phi$ is $(-)$twistorial. \qed
\end{thm} 

\indent 
Another consequence of Theorem \ref{thm:4to3} is that a one-dimensional conformal foliation   
with nowhere degenerate leaves on an oriented four-dimensional Riemannnian manifold is 
both $(+)$twistorial and $(-)$twistorial (Remark \ref{rem:tm}(5)\,) if and only if it is locally
generated by conformal vector fields \cite{Cal-sds} (the `if' part follows also from \cite{JonTod}\,). 
In particular, (the complexification of) 
any one-dimensional homothetic foliation locally defined by harmonic morphisms on an oriented  
four-dimensional Riemannian manifold is both $(+)$twistorial and $(-)$twistorial.\\  
\indent
With the same notations as above and in Section \ref{section:facts}, we have the following  
consequences of Theorem \ref{thm:4to3}, (in which we do not claim that (i)$\iff$(iii) is new). 

\begin{cor} \label{cor:4to3Hintegrable} 
Let $\V$ be a one-dimensional conformal foliation with integrable orthogonal complement 
on an oriented four-dimensional conformal manifold $(M^4,c)$\,.\\ 
\indent 
Then the following assertions are equivalent:\\ 
\indent 
{\rm (i)} $\V$ is $(\pm)$twistorial;\\ 
\indent 
{\rm (ii)} $\V$ produces harmonic morphisms with respect to any local representative 
of\/ $c$ for which $\V$ has geodesic leaves;\\
\indent 
{\rm (iii)} locally, there exist representatives of\/ $c$ with respect to which $\V$ is defined 
by totally geodesic Riemannian submersions. \qed 
\end{cor} 

\indent
We end this section with a result on $(\pm)$twistorial maps from conformally flat four-manifolds. 
Recall (see \cite{WarWel}) that, any conformally flat four-manifold can be locally embedded in the 
four-dimensional nonsingular quadric $Q_4$\,. 

\begin{prop} 
Let $M^4$ be an oriented conformally flat four-manifold and let $\phi:M^4\to N^3$ be a horizontally 
conformal submersion with nowhere degenerate fibres. Then, the following assertions are equivalent:\\ 
\indent
{\rm (i)} $\phi$ is $(-)$twistorial;\\ 
\indent
{\rm (ii)} $\phi$ is $(+)$twistorial.\\ 
\indent
Furthermore, if $\phi$ is $(\pm)$twistorial and $M^4$ is locally embedded in $Q_4$\,, then the twistor 
space $Z(N)$ of the induced Einstein--Weyl connection on $N^3$ is, locally, a surface in $Q_4$ such 
that the space of horizontal null geodesics of\/ $\phi$ is equal to the space of null geodesics on $Q_4$ 
which pass through $Z(N)\subseteq Q_4$\,. 
\end{prop}

\section{The classification of twistorial harmonic morphisms with one-dimensional fibres from self-dual four-manifolds} \label{section:4to3harmorph}  
 
\indent
In the last two sections we shall work in the real-analytic category.\\ 
\indent
Firstly, we list the four types of $(-)$twistorial harmonic morphisms with one-dimensional fibres which 
can be defined on a four-dimensional Riemannian manifold. Later on we shall see how they come from equation  
\eqref{e:m}\,.\\ 

{\bf Type 1 (`Killing type' \cite{Bry}).} Harmonic morphisms $\phi:M^4\to N^3$ whose fibres are
locally generated by nowhere zero Killing vector fields (see Section \ref{section:facts} above).\\ 

{\bf Type 2 (`warped-product type' \cite{BaiEel}).} Horizontally homothetic submersions with 
geodesic fibres orthogonal to an umbilical foliation by hypersurfaces (see Remark \ref{rem:1.9}(1)\,).\\

{\bf Type 3 (}cf.\ {\bf \cite{Pan-thesis,Pan-4to3,PanWoo-exm}).} $\phi:(M^4,g)\to(N^3,h)$ is a real-analytic
map which is, locally,
the canonical projection $\R\times N^3\to N^3$\,, $(\rho,x)\mapsto x$\,, and  
\begin{equation} \label{e:type3} 
g=\rho\,h+\rho^{-1}(\dif\!\rho+A)^2 
\end{equation} 
where $A$ is a one-form on $N^3$ which satisfies the Beltrami fields equation\linebreak $\dif\!A=-*A$.\\ 

{\bf Type 4 (}cf.\ \cite{Cal-sds}{\bf ).} $(N^3,h)$ is endowed with a Weyl connection $D$\,, $\phi:(M^4,g)\to(N^3,h)$
is a real-analytic map which is, locally, the canonical projection
$\R\times N^3\to N^3$\,, $(\rho,x)\mapsto x$\,, and  
\begin{equation} \label{e:type4} 
g=(e^{\rho}+c)\,h+\frac{1}{e^{\rho}+c}\,(\dif\!\rho-\a)^2 
\end{equation} 
where the Lee form $\a$ of $D$ with respect to $h$ satisfies the equation 
$$\dif\!\a-c*\a+*\dif\!c=0$$ on $N^3$ with $c:N^3\to\R$ a function (if $e^{\rho}+c<0$\,, 
then we replace $g$ with $-g$).\\ 

\begin{rem} \label{rem:warped-prod}
1) Note that maps of types 3 and 4 are always harmonic morphisms by Theorem \ref{thm:Bry}\,.\\
\indent
2) Let $\phi:(M^{n+1},g)\to(N^n,h)$\ $(n\geq1)$\ be a submersive harmonic morphism with
one-dimensional fibres. Then, the components of the fibres of $\phi$ form a homothetic foliation 
if and only if either $\phi$ is of type 1 or, up to a conformal change of $(M^{n+1},g)$ with factor 
constant along the fibres, $\phi$ is of type 2.\\  
\indent
3) If $\phi:(M^4,g)\to(N^3,h)$ is a harmonic morphism of type 3 then $\rho$ and $\phi^*(A)$ are
globally defined on $M^4$. Indeed, if we denote by $\l$ the dilation of $\phi$ then $\rho=\l^{-2}$ 
and $\phi^*(A)=-\dif\!^{\H}\bigr(\l^{-2}\bigl)$\,.\\ 
\indent
4) If $\phi$ is a harmonic morphism of type 4 with $\dif\!\a=0$ then the horizontal distribution
is integrable and, after a conformal change with basic factor, $\phi$ is of type 2\,. 
If $\a\neq0$\,, then $\rho$ and $c\circ\phi$ are globally defined on $M^4$ 
(this follows from the fact that, up to signs, we have 
$\theta=\dif\!\rho-\a$ and $\l^{-2}=e^{\rho}+c\circ\phi$ where $\l$ is the dilation of $\phi$); 
in particular, $e^{\rho}+c\circ\phi$ has constant sign on $M^4$.\\ 
\indent 
In the complex-analytic category the construction of type 4 is globally defined if $e^{\rho}+c\circ\phi$ is 
nowhere zero. Then, if $\a\neq0$\,, $e^{\rho}$ and $c\circ\phi$ are well-defined functions on $M^4$\,.
\end{rem}    

\begin{prop} \label{prop:c=pm1}
Let $\phi:(M^4,g)\to(N^3,h,D)$ be a harmonic morphism of type 4 with $c\neq0$\,. 
Then, outside the zero set of\/ $c$\,, after a conformal change with factor constant along the 
fibres and a suitable change of coordinates we may assume that $c=\pm1$ --- that is, $\phi$ is, 
locally, the canonical projection $\R\times N^3\to N^3$\,, $(\rho,x)\mapsto x$\,, and  
\begin{equation} \label{e:type4Beltrami} 
g=(e^{\rho}\pm1)\,h+\frac{1}{e^{\rho}\pm1}\,(\dif\!\rho-\a)^2 
\end{equation} 
where the Lee form $\a$ of\/ $D$ with respect to $h$ satisfies the Beltrami fields equation 
\begin{equation} \label{e:c=pm1} 
\dif\!\a=\pm*\a 
\end{equation} 
on $N^3$. 
\end{prop} 
\begin{proof} 
We shall show that $\widetilde{g}=|c|\,g$ is as in \eqref{e:type4Beltrami}\,, \eqref{e:c=pm1} after 
the change of coordinates $\widetilde{\rho}=\rho-\log|c|$\,.\\ 
\indent
Assume for simplicity that $c>0$ and let $\widetilde{h}=c^2h$. Then, the dilation $\widetilde{\l}$, 
fundamental vector field $\widetilde{V}$ and connection form $\widetilde{\theta}$ of the harmonic morphism 
$\phi:(M^4,\widetilde{g})\to(N^3,\widetilde{h})$ are as follows: $\widetilde{\l}^2=c\l^2$, $\widetilde{V}=V$ and $\widetilde{\theta}=\theta$. Also the 
Lee form of $D$ with respect to $\widetilde{h}$ is $\widetilde{\a}=\a-\dif\log c$.\\ 
\indent
It is easy to see that $$\dif\!\a-c*\a+*\dif\!c=0$$ written with respect to $\widetilde{h}$ becomes
$*\dif\!\widetilde{\a}=c\,\widetilde{\a}$. But, if we denote by $\widetilde{*}$ the Hodge-star operator on 
$(N^3,\widetilde{h})$, then $\widetilde{*}|_{\Lambda^2}=c^{-1}*\!|_{\Lambda^2}$ and hence 
\begin{equation} \label{e:tildeBeltrami} 
\widetilde{*}\dif\!\widetilde{\a}=\widetilde{\a}\;. 
\end{equation} 
\indent
Now it is easy to see that, with respect to $\widetilde{\rho}$, $\widetilde{h}$ and 
$\widetilde{\a}$, the metric $\widetilde{g}$ has the required form. 
\end{proof}  

\indent
Now we are able to state the main result of this section. 

\begin{thm} \label{thm:maint} 
Let $M^4$ and $N^3$ be real-analytic Riemannian manifolds of dimensions four and three, respectively. 
Suppose that $N^3$ is endowed with a Weyl connection $D$.\\ 
\indent
Then $\phi:M^4\to(N^3,D)$ is a $(-)$twistorial harmonic morphism if and only if, up to conformal 
changes of the metric on the domain and codomain, $\phi$ is of type 1\,,\,2\,,\,3 or 4\,; further,   
the conformal factor on $M^4$ can be taken to be constant along the fibres of\/ $\phi$\,.  
If $\phi$ is of type 2\,, then $D$ is the Levi-Civita connection of some metric in the conformal class of 
$N^3$, whilst if $\phi$ is of type 3\,, then $D$ is the Levi-Civita connection of\/ $N^3$. 
\end{thm} 
\begin{proof} 
By Corollary \ref{cor:4to3harmorph}, equation \eqref{e:m} is satisfied. If $V(\l^{-2})=0$ then 
$\phi$ is of type 1\,. {}From \cite[Corollary 1.5]{PanWood-d} it follows that $V(\l^{-2})$ 
is real-analytic,   
and from now on we shall assume that $V(\l^{-2})$ is nowhere zero on some open set. Moreover, by replacing, 
if necessary, $V$ with $-V$, we can assume that $V(\l^{-2})$ is positive at each point.\\ 
\indent
Let $X$ be a basic vector field. Note that the left-hand side of \eqref{e:m} must be basic and hence, 
using $[V,X]=0$ we obtain  
$X(V(\l^{-2}))-V(\l^{-2})\a(X)=0$\,. Thus 
\begin{equation} \label{e:Lee} 
\a=\dif\!^{\H}(\log V(\l^{-2}))\;. 
\end{equation} 
This implies that $X(\log V(\l^{-2}))$ is basic. Hence  
$$X(V(\log V(\l^{-2})))=V(X(\log V(\l^{-2})))=0\;.$$ 
\indent 
It follows that, if $V(\log V(\l^{-2}))$ is nonconstant then $\H$ is integrable, equivalently $\O=0$ on an open set, 
and hence, by real-analyticity, on $M^4$.  
This, together with \eqref{e:m}, gives $\a=-2\dif\!^{\H}(\log\l)$\,; hence $\dif\!^{\H}(\log\l)$ is basic,  
that is, $V(X(\log\l))=0$ for any basic vector field $X$. Because $\H$ is integrable this implies that $\V$ 
is homothetic, and hence after a conformal change with basic factor we get that $\phi$ is of type 2 
(see Remark \ref{rem:warped-prod}(2)\,).
Moreover, as $\a$ is exact, $D$ is the Levi-Civita connection of some Riemannian metric 
in the conformal class of $N^3$.\\    
\indent
{}From now on we shall suppose that $V(\log V(\l^{-2}))=a$ for some constant $a\in\R$\,. 
If $a=0$ this implies that $V(\l^{-2})$ is basic. We can write $V=\partial/\partial\!\rho$ for some 
function $\rho$ which is zero on a chosen section of $\phi$\,. Hence $\l^{-2}=b\rho+c$ for some basic 
functions $b$ and $c$\,. 
After suitable conformal changes on $N^3$ and $M^4$ we have $\l^{-2}=\rho+c/b$\,, i.e., 
$V(\l^{-2})=1$\,. By \eqref{e:Lee}\,, this implies that $\a=0$ and hence $D$ is the Levi-Civita 
connection of $N^3$. 
Moreover, we can locally write $\theta=\dif(\l^{-2})+\phi^*(A)$ for some one-form $A$ on $N^3$. 
Hence $\dif\!^{\H}(\l^{-2})=-\phi^*(A)$ and $\O=\dif\!\theta=\phi^*(\dif\!A)$ which together 
with \eqref{e:m} gives $-A=*\dif\!A$\,. Thus we have proved that, if $a=0$ then, up to a
conformal change with basic factor, $\phi$ is of type 3\,.\\ 
\indent
It remains to consider the case when $V(\log V(\l^{-2}))=a$ for some nonzero constant $a$\,. 
Then, if we set $\rho=\log V(\l^{-2})$\,, we can assume that $V=a\,\partial/\partial\!\rho$\,. 
As $V(\l^{-2})=e^{\rho}$ we have that $a\l^{-2}=e^{\rho}+c$ for some basic function $c$\,. 
Hence 
\begin{equation} 
\l^{-2}=\frac{1}{a}\,(e^{\rho}+c)\;. 
\end{equation} 
\indent
{}From $\dif\!^{\V}\rho=a\theta$ and \eqref{e:Lee} we get 
$\dif\!\rho=a\theta+\a$\,; equivalently 
\begin{equation} \label{e:theta} 
\theta=\frac{1}{a}\,\bigl(\dif\!\rho-\a\bigr)\;. 
\end{equation}
\indent
Now, we can write 
\begin{equation} \label{e:nowwecanwrite} 
\begin{split} 
&\:a(\dif\!^{\H}-\a)(\l^{-2})=a\bigl(\dif\!^{\H}-\dif\!^{\H}\rho\bigr)\bigl(\frac1a\,(e^{\rho}+c)\bigr)\\ 
&=e^{\rho}\dif\!^{\H}\rho+\dif\!^{\H}c-e^{\rho}\dif\!^{\H}\rho-c\dif\!^{\H}\rho\\ 
&=\dif\!^{\H}c-c\dif\!^{\H}\rho=\dif\!c-c\dif\!^{\H}\rho=\dif\!c-c\a\;. 
\end{split} 
\end{equation} 
\indent
On the other hand, from \eqref{e:theta} it follows that $a\,\O=a\dif\!\theta=-\dif\!\a$ which,  together with \eqref{e:m}
and \eqref{e:nowwecanwrite}\,, gives $\dif\!c-c\a=-*\dif\!\a$\,.
Thus we have proved that, if $\phi$ is not of type 1\,,\,2 or 3 then 
\begin{equation} \label{e:alphacnonconst} 
\dif\!\a-c*\a+*\dif\!c=0 
\end{equation} 
and, locally, the metric of $M^4$ can be written in the form 
\begin{equation} \label{e:gcnonconst} 
g=\frac1a\,\bigl\{\,(e^{\rho}+c)\,h+\frac{1}{e^{\rho}+c}\,(\dif\!\rho-\a)^2\,\bigr\} 
\end{equation} 
where, obviously, we can assume that $a=\pm1$ and that $c$ is nonzero; since, if $c=0$, then 
after a suitable conformal change of basic factor, $\phi$ is of type 2\,. 
\end{proof} 

\indent
The following result is an immediate consequence of Theorem \ref{thm:4to3} and Theorem \ref{thm:maint}.  

\begin{cor} \label{cor:maintsd} 
Let $M^4$ be a self-dual manifold with a real-analytic metric and let $\phi:M^4\to(N^3,D)$ 
be a $(-)$twistorial harmonic morphism.\\ 
\indent
Then $D$ is Einstein--Weyl and, up to conformal changes of the metrics, 
$\phi$ is of type 1\,,\,2\,,\,3 or 4\,; further,   
the conformal factor on $M^4$ can be taken to be constant along the fibres of\/ $\phi$\,.  
If $\phi$ is of type 2\,, then both 
$M^4$ and $N^3$ are conformally flat, whilst if $\phi$ is of type 3\,, then $N^3$ has constant curvature. 
\end{cor}

\section{Constructions of self-dual metrics}  \label{section:newconstr}

\indent
In this section we show that harmonic morphisms are related to known constructions
of self-dual metrics.\\ 

\indent
Firstly, we recall the following construction of P.E.~Jones and K.P.~Tod \cite{JonTod} 
(see \cite{LeB-jdg}\,, cf.\ \cite{GibHaw,Haw}). 

\begin{thm} \label{thm:JonTodconstr} 
Let $(N^3,[h],D)$ be an Einstein--Weyl 3-manifold and let $(M^4,N^3,S^1)$ be a (local) principal bundle  
endowed with a (local) principal connection $\H\subseteq TM$\,. Define a Riemannian metric 
$g$ on $M^4$ by 
\begin{equation} 
g=v\,\phi^*(h)+v^{-1}\,\theta^2 
\end{equation} 
where $\phi:M^4\to N^3$ is the projection of the principal bundle $(M^4,N^3,S^1)$\,, $v=u\circ\phi$ for 
some positive smooth function $u$ on $N^3$ and $\theta$ is the connection form of\/ $\H$\,.\\ 
\indent
Then $(M^4,g)$ is self-dual (respectively, anti-self-dual) 
if the following $S^1$-monopole equation holds on $N^3$: 
\begin{equation} 
(\dif-\a)\,u=*F\qquad\textit{(respectively,\:\,$(\dif-\a)\,u=-*F$\,)}    
\end{equation} 
where $\a$ is the Lee form of\/ $D$ with respect to $h$ and $F\in\G(\Lambda^2(T^*N))$ is the 
curvature form of\/ $\H$\,. 
\end{thm} 

\begin{rem} 
The construction of Theorem \ref{thm:JonTodconstr} clearly gives harmonic morphisms of type 1 and, if 
$\H$ is flat, of type 2, up to a conformal change with basic factor. 
\end{rem} 
 
\begin{thm}[{\bf Type 3}, cf.\ \cite{Cal-sds}, \cite{PanWoo-exm}] \label{thm:3sd}
Let\/ $(N^3,h)$ be a constant curvature 3-mani\-fold and let\/ $A$ be a 
one-form on $N^3$. Define a Riemannian metric on\/ $(0,\infty)\times N^3$ by 
\begin{equation} \label{e:3sdconstruction} 
g=\rho\,h+\rho^{-1}(\dif\!\rho+A)^2\qquad\qquad(\rho\in(0,\infty)\,)\;. 
\end{equation} 
\indent
Then\/ $g$ is self-dual (respectively, anti-self-dual) if the following 
Beltrami 
fields equation holds on\/ $N^3$: 
\begin{equation} \label{e:3sd} 
\dif\!A=-*A\qquad\textit{(respectively,\:\,$\dif\!A=*A$\,)\;.} 
\end{equation} 
\end{thm} 

\begin{thm}[{\bf Type 4}, \cite{Cal-sds}] \label{thm:4sd}
Let $(N^3,[h],D)$ be an Einstein--Weyl 3-manifold. Define a Riemannian  
metric on\/ $(0,\infty)\times N^3$ by 
\begin{equation}
g=(e^{\rho}\pm1)\,h+\frac{1}{e^{\rho}\pm1}\,(\dif\!\rho-\a)^2\qquad\qquad(\rho\in(0,\infty)\,)
\end{equation} 
where $\a$ is the Lee form of\/ $D$ with respect to $h$\,.\\ 
\indent
Then $g$ is self-dual (respectively, anti-self-dual) if the following Beltrami 
fields equation holds on $N^3$ 
\begin{equation} \label{e:4sd}  
\dif\!\a=\pm*\a\qquad\textit{(respectively,\:\,$\dif\!\a=\mp*\a$\,)}\;.
\end{equation}  
\end{thm} 

\begin{rem}
1) The construction of Theorem \ref{thm:JonTodconstr} gives Einstein metrics if $\a=0$ and $h$ is flat,
in which case $g$ is Ricci-flat self-dual \cite{GibHaw,Haw} (see \cite{PanWoo-exm}).\\ 
\indent 
2) The construction of Theorem \ref{thm:3sd} gives Einstein metrics if and only if $(N^3,h)$ has constant 
sectional curvature equal to $1/4$, in which case $g$ is Ricci-flat self-dual. The construction of 
Theorem \ref{thm:4sd} gives Einstein metrics if and only if $\a=0$ and so $h$ has constant 
curvature, in which case $g$ also has constant curvature. These facts follow from 
\cite[Theorem 1.5]{PanWoo-exm} and are also true if we make a conformal change of $g$ with 
factor constant along the fibres of $\phi$\,.\\ 
\indent 
3) Let $A$ be a solution of the Beltrami fields equation \eqref{e:3sd} on a three-dimension\-al Riemannian 
manifold $(N^3,h)$ with constant curvature. Then the components of $A$ with respect to an orthonormal 
basis of left invariant one-forms are eigenfunctions of the Laplace-Beltrami operator of $(N^3,h)$ 
(cf.\  \cite{PanWoo-exm}\,; note that the corresponding eigenvalues are imaginary if $(N^3,h)$ 
has negative sectional curvature).\\ 
\indent
4) Let $(N^3,[h],D)$ be an Einstein--Weyl space for which there exists a surjective submersion 
$\pi:Z(N)\to\C\!P^1$ whose fibres are transversal to the twistor lines. Then, it is well known 
that, for any $x_1\,,\,x_2\in\C\!P^1=S^2$\,, the angle formed by any leaf of the 
foliation corresponding to $\pi^{-1}(x_1)$ and any leaf of the foliation corresponding to $\pi^{-1}(x_2)$ 
is equal to $\mathrm{dist}_{S^2}(x_1,x_2)$ (see \cite{GauTod}\,, \cite{CalPed}\,).\\
\indent 
For example, if $D$ is the Levi-Civita connection of a Riemannian manifold $(N^3,h)$ with constant 
nonnegative sectional curvature then $(N^3,[h],D)$ has this property (see \cite{BaiWoo2}\,). 
It follows that, if the codomain of a harmonic morphism of type 3 has nonnegative constant sectional 
curvature then its domain is hyper-Hermitian.\\
\indent
5) D.M.J.~Calderbank \cite{Cal-sds} gives the type 4 construction with the extra condition on
$(N^3,[h],D)$ that it is an Einstein--Weyl space for which there exists a surjective submersion
$\pi:Z(N)\to\C\!P^1$ whose fibres are transversal to the twistor lines. Then the construction
gives hyper-Hermitian metrics which, after a suitable conformal change, are Einstein with
nonzero scalar curvature.
\end{rem}


\begin{thebibliography}{10} 

                             


\bibitem{AtHiSi} 
M.F.~Atiyah, N.J.~Hitchin, I.M.~Singer, Self-duality in four-dimensional 
Riemannian geometry, \textit{Proc. Roy. Soc. London Ser. A}, {\bf 362} 
(1978) 425--461. 
\bibitem{BaiEel}
P.~Baird, J.~Eells, A conservation law for harmonic maps,  
\textit{Geometry Symposium. Utrecht 1980}, 1--25, Lecture Notes in Math. 894, 
Springer-Verlag, Berlin, 1981. 
\bibitem{BaiWoo-Bernstein} 
P.~Baird, J.C.~Wood, Bernstein theorems for harmonic morphisms from $R\sp 3$ and $S\sp 3$,  
\textit{Math. Ann.}, {\bf 280} (1988) 579--603. 
\bibitem{BaiWoo-spfo}
P.~Baird, J.C.~Wood, Harmonic morphisms and conformal foliations by
geodesics of three-dimensional 
space forms, \emph{J. Austral. Math. Soc. (A)}, {\bf 51} (1991) 118--153. 
\bibitem{BaiWoo1}
P.~Baird, J.C.~Wood, Harmonic morphisms, Seifert fibre spaces and conformal 
foliations, \textit{Proc. London Math. Soc.}, {\bf 64} (1992) 170--196.  
\bibitem{BaiWoo2}
P.~Baird, J.C.~Wood, \textit{Harmonic morphisms between Riemannian manifolds}, 
London Math. Soc. Monogr. (N.S.), no. 29, Oxford Univ. Press., Oxford, 2003.  
\bibitem{Bott-partial} 
R.~Bott,  
Lectures on characteristic classes and foliations. 
Notes by Lawrence Conlon, with two appendices by J. Stasheff,  
\textit{Lectures on algebraic and differential topology (Second Latin American School in Math., Mexico City, 1971)}, 1--94. Lecture Notes in Math. 279, Springer, Berlin, 1972.  
\bibitem{Bry}
R.L.~Bryant, Harmonic morphisms with fibres of dimension one, \textit{Comm. 
Anal. Geom.}, {\bf 8} (2000) 219--265. 
\bibitem{Cal-sds} 
D.M.J.~Calderbank, Selfdual Einstein metrics and conformal submersions, 
\textit{Preprint}, Edinburgh University, 2000; ArXiv math DG/0001041. 
\bibitem{CalPed} 
D.M.J.~Calderbank, H.~Pedersen, Selfdual spaces with complex structures, Einstein--Weyl 
geometry and geodesics, \textit{Ann. Inst. Fourier (Grenoble)}, {\bf 50} (2000) 921--963. 
\bibitem{EelSal} 
J.~Eells, S.~Salamon, Twistorial construction of harmonic maps of surfaces into four-manifolds, \textit{Ann. Scuola Norm. Sup. Pisa Cl. Sci. (4)}, {\bf 12} (1985) 589--640. 
\bibitem{Fug}
B.~Fuglede, Harmonic morphisms between Riemannian manifolds, \textit{Ann. 
Inst. Fourier (Grenoble)}, {\bf 28} (1978) 107--144. 
\bibitem{Gau} 
P.~Gauduchon, Structures de Weyl-Einstein, espaces de twisteurs et vari\'et\'es de type $S^1\times S^3$, 
\textit{J. Reine Angew. Math.}, {\bf 469} (1995) 1--50. 
\bibitem{GauTod} 
P.~Gauduchon, K.P.~Tod, Hyper-Hermitian metrics with symmetry, \textit{J. Geom. Phys.}, {\bf 25} (1998) 291--304. 
\bibitem{GibHaw} 
G.W.~Gibbons, S.W.~Hawking, Gravitational multi-instantons, 
\textit{Phys. Lett. B}, {\bf 78} (1978) 430--432.  
\bibitem{Gudbib}
S.~Gudmundsson, \textit{The Bibliography of Harmonic Morphisms}, 
{\tt http://www.maths.lth.se/\\
matematiklu/personal/sigma/harmonic/bibliography.html} 
\bibitem{Haw}
S.W.~Hawking, Gravitational Instantons, \textit{Phys. Lett. A}, {\bf 60} (1977) 81--83. 
\bibitem{Hit-complexmfds}  
N.J.~Hitchin, Complex manifolds and Einstein's equations, 
\textit{Twistor geometry and nonlinear systems (Primorsko, 1980)}, 73--99, 
Lecture Notes in Math., 970, Springer, Berlin, 1982. 
\bibitem{Ish}
T.~Ishihara, A mapping of Riemannian manifolds which  preserves 
harmonic functions, \textit{J. Math. Kyoto Univ.}, {\bf 19} (1979) 215--229. 
\bibitem{JonTod} 
P.E.~Jones, K.P.~Tod, Minitwistor spaces and Einstein--Weyl spaces, 
\textit{Classical Quantum Gravity}, {\bf 2} (1985) 565--577. 
\bibitem{KenPlu} 
P.C.~Kendall, C.~Plumpton, \textit{Magnetohydrodynamics, with Hydrodynamics}, 
Pergamon Press, Oxford; Macmillan, 1964. 
\bibitem{KoNo}
S.~Kobayashi, K.~Nomizu, \textit{Foundations of differential geometry}, I, II, 
Wiley Classics Library (reprint of the 1963, 1969 original), Wiley-Interscience Publ.,  
Wiley, New-York, 1996. 
\bibitem{Kod} 
K.~ Kodaira, A theorem of completeness of characteristic systems for analytic families of 
compact submanifolds of complex manifolds, \textit{Ann. of Math.}, {\bf 75} (1962) 146--162. 
\bibitem{LeB-thesis} 
C.R.~LeBrun, Spaces of complex geodesics and related structures, Ph.D. Thesis, Oxford, 1980.  
\bibitem{LeB-nullgeod} 
C.R.~LeBrun, Spaces of complex null geodesics in complex-Riemannian geometry, 
\textit{Trans. Amer. Math. Soc.}, {\bf 278} (1983) 209--231.  
\bibitem{LeB-jdg} 
C.R.~LeBrun, Explicit self-dual metrics on $\C\!P^2\#\cdots\#\C\!P^2$, 
\textit{J. Differential Geom.}, {\bf 34} (1991) 223--253. 
\bibitem{TN-III} 
L.J.~Mason, L.P.~Hughston, P.Z.~Kobak, K.~Pulverer (Editors), 
\textit{Further advances in twistor theory, Vol. III: Curved twistor spaces},
Chapman \& Hall/CRC Research Notes in Mathematics, 424. Chapman \& Hall/CRC, Boca Raton, FL, 2001.
\bibitem{MasWoo} 
L.J.~Mason, N.M.J.~Woodhouse, \textit{Integrability, self-duality, and twistor theory}, 
London Math. Soc. Monogr. (N.S.), no. 15, Oxford Univ. Press., Oxford, 1996.     
\bibitem{Pan}
R.~Pantilie, Harmonic morphisms with one-dimensional fibres, 
\textit{Internat. J. Math.}, {\bf 10} (1999) 457--501. 
\bibitem{Pan-thesis} 
R.~Pantilie, Submersive harmonic maps and morphisms, Ph.D. Thesis, University of Leeds, 2000. 
\bibitem{Pan-4to3} 
R.~Pantilie, Harmonic morphisms with 1-dimensional fibres on 4-dimensional  Einstein manifolds, 
\textit{Comm. Anal. Geom.}, {\bf 10} (2002) 779--814.  
\bibitem{PanWood-d} 
R.~Pantilie, J.C.~Wood, Harmonic morphisms with one-dimensional fibres on Einstein manifolds, 
\textit{Trans. Amer. Math. Soc}, {\bf 354} (2002) 4229--4243. 
\bibitem{PanWoo-exm} 
R.~Pantilie, J.C.~Wood, A new construction of Einstein self-dual manifolds, 
\textit{Asian J. Math.}, {\bf 6} (2002) 337--348.
\bibitem{Pen} 
R.~Penrose, Nonlinear gravitons and curved twistor theory, 
\textit{Gen. Relativity Gravitation}, {\bf 7} (1976) 31--52. 
\bibitem{WarWel} 
R.S.~Ward, R.O.~Wells, \textit{Twistor geometry and field theory}, Cambridge University Press, 
Cambridge, 1990. 
\bibitem{Woo}
J.C.~Wood, Harmonic morphisms, foliations and Gauss maps, \textit{Complex 
differential geometry and nonlinear differential equations}, 145--183, 
Contemp. Math. 49, Amer Math. Soc., Providence, RI, 1986. 
\bibitem{Woo-4d} 
J.C.~Wood, Harmonic morphisms and Hermitian structures on Einstein $4$-manifolds. 
\textit{Internat. J. Math.}, {\bf 3} (1992) 415--439. 
\end{thebibliography}
\end{document}